\documentclass[final,preprint,3p,sort&compress]{elsarticle}
\journal{{Journal of Computational Physics}}

\usepackage[dvips]{epsfig}
\usepackage{graphicx} 
\usepackage{pdfpages}
\usepackage{latexsym}
\usepackage{verbatim}
\usepackage{amsmath,amsthm}
\usepackage{amssymb}
\usepackage{bbm}
\usepackage{fonts}
\usepackage{enumerate}
\usepackage{placeins}
\usepackage{standalone}
\usepackage{paralist}
\usepackage{subcaption}
\definecolor{greenyellow}   {cmyk}{0.15, 0   , 0.69, 0   }
\definecolor{yellow}        {cmyk}{0   , 0   , 1   , 0   }
\definecolor{goldenrod}     {cmyk}{0   , 0.10, 0.84, 0   }
\definecolor{dandelion}     {cmyk}{0   , 0.29, 0.84, 0   }
\definecolor{apricot}       {cmyk}{0   , 0.32, 0.52, 0   }
\definecolor{peach}         {cmyk}{0   , 0.50, 0.70, 0   }
\definecolor{melon}         {cmyk}{0   , 0.46, 0.50, 0   }
\definecolor{yelloworange}  {cmyk}{0   , 0.42, 1   , 0   }
\definecolor{orange}        {cmyk}{0   , 0.61, 0.87, 0   }
\definecolor{burntorange}   {cmyk}{0   , 0.51, 1   , 0   }
\definecolor{bittersweet}   {cmyk}{0   , 0.75, 1   , 0.24}
\definecolor{redorange}     {cmyk}{0   , 0.77, 0.87, 0   }
\definecolor{mahogany}      {cmyk}{0   , 0.85, 0.87, 0.35}
\definecolor{maroon}        {cmyk}{0   , 0.87, 0.68, 0.32}
\definecolor{brickred}      {cmyk}{0   , 0.89, 0.94, 0.28}
\definecolor{red}           {cmyk}{0   , 1   , 1   , 0   }
\definecolor{orangered}     {cmyk}{0   , 1   , 0.50, 0   }
\definecolor{rubinered}     {cmyk}{0   , 1   , 0.13, 0   }
\definecolor{wildstrawberry}{cmyk}{0   , 0.96, 0.39, 0   }
\definecolor{salmon}        {cmyk}{0   , 0.53, 0.38, 0   }
\definecolor{carnationpink} {cmyk}{0   , 0.63, 0   , 0   }
\definecolor{magenta}       {cmyk}{0   , 1   , 0   , 0   }
\definecolor{violetred}     {cmyk}{0   , 0.81, 0   , 0   }
\definecolor{rhodamine}     {cmyk}{0   , 0.82, 0   , 0   }
\definecolor{mulberry}      {cmyk}{0.34, 0.90, 0   , 0.02}
\definecolor{redviolet}     {cmyk}{0.07, 0.90, 0   , 0.34}
\definecolor{fuchsia}       {cmyk}{0.47, 0.91, 0   , 0.08}
\definecolor{lavender}      {cmyk}{0   , 0.48, 0   , 0   }
\definecolor{thistle}       {cmyk}{0.12, 0.59, 0   , 0   }
\definecolor{orchid}        {cmyk}{0.32, 0.64, 0   , 0   }
\definecolor{darkorchid}    {cmyk}{0.40, 0.80, 0.20, 0   }
\definecolor{purple}        {cmyk}{0.45, 0.86, 0   , 0   }
\definecolor{plum}          {cmyk}{0.50, 1   , 0   , 0   }
\definecolor{violet}        {cmyk}{0.79, 0.88, 0   , 0   }
\definecolor{royalpurple}   {cmyk}{0.75, 0.90, 0   , 0   }
\definecolor{blueviolet}    {cmyk}{0.86, 0.91, 0   , 0.04}
\definecolor{periwinkle}    {cmyk}{0.57, 0.55, 0   , 0   }
\definecolor{cadetblue}     {cmyk}{0.62, 0.57, 0.23, 0   }
\definecolor{cornflowerblue}{cmyk}{0.65, 0.13, 0   , 0   }
\definecolor{midnightblue}  {cmyk}{0.98, 0.13, 0   , 0.43}
\definecolor{navyblue}      {cmyk}{0.94, 0.54, 0   , 0   }
\definecolor{royalblue}     {cmyk}{1   , 0.50, 0   , 0   }
\definecolor{blue}          {cmyk}{1   , 1   , 0   , 0   }
\definecolor{cerulean}      {cmyk}{0.94, 0.11, 0   , 0   }
\definecolor{cyan}          {cmyk}{1   , 0   , 0   , 0   }
\definecolor{processblue}   {cmyk}{0.96, 0   , 0   , 0   }
\definecolor{skyblue}       {cmyk}{0.62, 0   , 0.12, 0   }
\definecolor{turquoise}     {cmyk}{0.85, 0   , 0.20, 0   }
\definecolor{tealblue}      {cmyk}{0.86, 0   , 0.34, 0.02}
\definecolor{aquamarine}    {cmyk}{0.82, 0   , 0.30, 0   }
\definecolor{bluegreen}     {cmyk}{0.85, 0   , 0.33, 0   }
\definecolor{emerald}       {cmyk}{1   , 0   , 0.50, 0   }
\definecolor{junglegreen}   {cmyk}{0.99, 0   , 0.52, 0   }
\definecolor{seagreen}      {cmyk}{0.69, 0   , 0.50, 0   }
\definecolor{green}         {cmyk}{1   , 0   , 1   , 0   }
\definecolor{forestgreen}   {cmyk}{0.91, 0   , 0.88, 0.12}
\definecolor{pinegreen}     {cmyk}{0.92, 0   , 0.59, 0.25}
\definecolor{limegreen}     {cmyk}{0.50, 0   , 1   , 0   }
\definecolor{yellowgreen}   {cmyk}{0.44, 0   , 0.74, 0   }
\definecolor{springgreen}   {cmyk}{0.26, 0   , 0.76, 0   }
\definecolor{olivegreen}    {cmyk}{0.64, 0   , 0.95, 0.40}
\definecolor{rawsienna}     {cmyk}{0   , 0.72, 1   , 0.45}
\definecolor{sepia}         {cmyk}{0   , 0.83, 1   , 0.70}
\definecolor{brown}         {cmyk}{0   , 0.81, 1   , 0.60}
\definecolor{tan}           {cmyk}{0.14, 0.42, 0.56, 0   }
\definecolor{gray}          {cmyk}{0   , 0   , 0   , 0.50}
\definecolor{black}         {cmyk}{0   , 0   , 0   , 1   }
\definecolor{white}         {cmyk}{0   , 0   , 0   , 0   } 

 \usepackage[]{hyperref} 
        
\usepackage{pgfplots}
\pgfplotsset{compat=newest}       
\newcommand{\externaltikz}[2]{\includegraphics{Externals/#1}}		

 \usepackage{relinput}
\newtheorem{theorem}{Theorem}[section]
\newtheorem{definition}[theorem]{Definition}
\newtheorem{remark}[theorem]{Remark}
\newtheorem{example}[theorem]{Example}
\newtheorem{assumption}[theorem]{Assumption}

\newtheorem{lemma}[theorem]{Lemma}
\newtheorem{corollary}[theorem]{Corollary}

\newcounter{tikzsubfigcounter}[figure]
\renewcommand{\thetikzsubfigcounter}{\the\numexpr\value{figure}+1\relax\alph{tikzsubfigcounter}}

\newcounter{tikzsubfigcounterinvisible}[figure]
\renewcommand{\thetikzsubfigcounterinvisible}{\the\numexpr\value{figure}+1\relax\alph{tikzsubfigcounterinvisible}}

\newcommand{\settikzlabel}[1]{ %
\refstepcounter{tikzsubfigcounterinvisible} \label{#1} 
}

\numberwithin{equation}{section}
\newcommand{\bdm}{\begin{displaymath}}
\newcommand{\edm}{\end{displaymath}}
\newcommand{\beq}{\begin{equation}}
\newcommand{\eeq}{\end{equation}}
\newcommand{\beqa}{\begin{eqnarray}}
\newcommand{\eeqa}{\end{eqnarray}}

\title{Kershaw closures for linear transport equations in slab geometry I: model derivation }
\author[fs]{Florian Schneider}
\address[fs]{Fachbereich Mathematik, TU Kaiserslautern, Erwin-Schr\"odinger-Str., 67663 Kaiserslautern, Germany, {\tt schneider@mathematik.uni-kl.de}}

\date{}

\usepgfplotslibrary{groupplots}

\newlength{\figureheight}
\newlength{\figurewidth}
\usetikzlibrary{arrows,shapes,positioning}
\usetikzlibrary{decorations.markings}
\tikzstyle arrowstyle=[scale=1]
\tikzstyle directed=[postaction={decorate,decoration={markings,
		mark=at position .65 with {\arrow[arrowstyle]{stealth}}}}]
\tikzstyle reverse directed=[postaction={decorate,decoration={markings,
		mark=at position .65 with {\arrowreversed[arrowstyle]{stealth};}}}]

\newcommand{\secref}[1]{Section~\ref{#1}}

\newcommand{\lemref}[1]{Lemma~\ref{#1}}

\newcommand{\rmref}[1]{Remark~\ref{#1}}
\renewcommand{\corref}[1]{Corollary~\ref{#1}}
\newcommand{\figref}[1]{Figure~\ref{#1}}

\newcommand{\abs}[1]{\ensuremath{\left| #1 \right|}}

\newcommand{\R}{\mathbb{R}}
\newcommand{\Rpos}{\R_{\geq 0}}
\newcommand{\scattering}{\ensuremath{\sigma_s}}
\newcommand{\absorption}{\ensuremath{\sigma_a}}

\newcommand{\source}{\ensuremath{Q}}

\newcommand{\sphere}[1][2]{\ensuremath{\mathcal{S}^{#1}}}

\newcommand{\distribution}[1][ ]{\ensuremath{\psi_{#1}}}
\newcommand{\distributiontzero}{\ensuremath{\distribution[\timevar=0]}}
\newcommand{\distributionboundary}{\ensuremath{\distribution[b]}}
\newcommand{\distributionvacuum}{\ensuremath{\distribution[\text{vac}]}}
\newcommand{\ansatz}[1][ ]{\ensuremath{\hat{\psi}_{#1}}}

\newcommand{\momentorder}{\ensuremath{N}}
\newcommand{\momentnumber}{\ensuremath{n}}
\newcommand{\basis}[1][ ]{{\ensuremath{\bb_{#1}}}} %Basis
\newcommand{\basisind}{\ensuremath{i}} %Basis
\newcommand{\basiscomp}[1][\basisind]{\ensuremath{b_{#1}}} %Basis
\newcommand{\normalizedbasis}[1][ ]{{\ensuremath{\widehat{\bb}_{#1}}}} %Basis
\newcommand{\fmbasis}[1][\momentorder]{\basis[#1]} %Basis
\newcommand{\moments}[1][ ]{\ensuremath{\bu_{#1}}} %moment vector
\newcommand{\momentcomp}[1]{\ensuremath{u_{#1}}} %moment vector

\newcommand{\isotropicmoment}{\moments[\text{iso}]}
\newcommand{\normalizedmoments}[1][ ]{\ensuremath{\bsphi_{#1}}} %moment vector
\newcommand{\normalizedmomentcomp}[1]{\ensuremath{\phi_{#1}}} %moment vector
\newcommand{\normalizedisotropicmoment}{\normalizedmoments[\text{iso}]}
\newcommand{\multipliers}[1][ ]{\ensuremath{\bsalpha_{#1}}} %moment vector
\newcommand{\convexscalar}{\ensuremath{\zeta}} %scalar indicating a convex combination
\newcommand{\flow}{\ensuremath{f_{\text{low}}}} %eddington factor
\newcommand{\fup}{\ensuremath{f_{\text{up}}}} %eddington factor

\newcommand{\SC}{\ensuremath{\Omega}} %moment vector
\newcommand{\SCheight}{\ensuremath{\mu}} %moment vector
\newcommand{\Domain}{\ensuremath{X}} %moment vector
\newcommand{\spatialVariable}{\ensuremath{\bx}} %moment vector
\newcommand{\timeint}{\ensuremath{T}} %time interval
\newcommand{\tf}{\ensuremath{t_f}} %final time
\newcommand{\timevar}{\ensuremath{t}} %final time

\newcommand{\ints}[1]{\ensuremath{\left<#1\right>}}

\newcommand{\collisionop}{\ensuremath{\cC}}
\newcommand{\collision}[1]{\ensuremath{\collisionop\left(#1\right)}}
\newcommand{\lincollisionop}{\ensuremath{\cI}}
\newcommand{\lincollision}[1]{\ensuremath{\lincollisionop\left(#1\right)}}
\newcommand{\collisionkernel}{\ensuremath{K}}

\newcommand{\dirac}{\ensuremath{\delta}}

\newcommand{\Lp}[1]{\ensuremath{L_{#1}}}
\newcommand{\RD}[2]{\ensuremath{\mathcal{R}_{#1}^{#2}}}
\newcommand{\RDone}[1]{\left. \RD{#1}{} \right|_{\momentcomp{0} = 1}}

\newcommand{\PN}[1][\momentorder]{\ensuremath{\text{P}_{#1}}}
\newcommand{\MN}[1][\momentorder]{\ensuremath{\text{M}_{#1}}}
\newcommand{\KN}[1][\momentorder]{\ensuremath{\text{K}_{#1}}}

\newcommand{\Flux}{\ensuremath{\bF}}
\newcommand{\Source}{\ensuremath{\bs}}

\newcommand{\eigenvector}{\ensuremath{\bv}}

\newcommand{\eigenvalue}{\ensuremath{\lambda}}

\DeclareMathOperator*{\argmin}{argmin}
\newcommand{\ld}[1]{\ensuremath{{#1}_*}} %Legendre dual
\newcommand{\entropy}{\ensuremath{\eta}} %Entropy function
\newcommand{\entropyFunctional}{\ensuremath{\mathcal{H}}} %Entropy functional
\def\quand{\quad \mbox{and} \quad}

\newcommand{\z}{\ensuremath{z}}

\newcommand{\dz}{\partial_{\z}}
\newcommand{\dt}{\partial_\timevar}

\newcommand{\hankelA}{A}
\newcommand{\hankelB}{B}
\newcommand{\hankelC}{C}

\newcommand{\hankelv}{v}
\newcommand{\hankelb}{\bsbeta} %Eventuell austauschen.........................................................
\newcommand{\hankelw}{\bw}
\newcommand{\hankelrho}{\rho}
\newcommand{\hankelphi}{\gamma}
\newcommand{\hankelrank}{r}
\newcommand{\hankelind}{i}
\newcommand{\hankelhalfind}{k}
\newcommand{\generatingFunction}[1][\hankelphi]{\ensuremath{g_{#1}}}
\newcommand{\range}[1]{R \left( #1 \right)}

\newcommand{\pseudoinv}[1]{\ensuremath{#1^\dag}}

\newcommand{\zL}{\ensuremath{\z_{L}}}
\newcommand{\zR}{\ensuremath{\z_{R}}}

\newcommand{\ncells}{\ensuremath{n_{\z}}}

\pgfplotscreateplotcyclelist{color parula}{% 
	royalblue!70,every mark/.append style={solid,line width = 0pt,fill=royalblue!60!black},mark=ball\\%
	goldenrod!50!yelloworange,every mark/.append style={solid,fill=goldenrod!30!black},mark=*\\%
	green,every mark/.append style={black,fill=yellowgreen},mark=diamond*\\%
	bluegreen,every mark/.append style={fill=black},mark=triangle*\\%	
	royalblue!70,densely dashed,every mark/.append style={solid,line width = 0pt,fill=royalblue!60!black},mark=ball\\%
	goldenrod!50!yelloworange,densely dashed,every mark/.append style={solid,black,fill=goldenrod},mark=diamond*\\%
	yellowgreen,densely dashed,every mark/.append style={solid,fill=yellowgreen!60!royalblue},mark=square*, mark size= 1pt\\%
	bluegreen,densely dashed,mark size = 1.5pt,every mark/.append style={solid,fill=white},mark=otimes*\\%,densely dashed	
	royalblue,densely dashed,mark size = 1.5pt,every mark/.append style={solid,fill=royalblue!30!black},mark=pentagon*\\%,densely dashed
	goldenrod!40!yelloworange,every mark/.append style={solid,fill=goldenrod!30!black},mark=|\\%		
	}

\begin{document}

\begin{abstract}
This paper provides a new class of moment models for linear kinetic equations in slab geometry. These models can be evaluated cheaply while preserving the important realizability property, that is the fact that the underlying closure is non-negative. Several comparisons with the (expensive) 	state-of-the-art minimum-entropy models are made, showing the similarity in approximation quality of the two classes. 
\end{abstract}
\maketitle

\noindent

% {\bf Key words.}

\section{Introduction}
In recent years, many approaches have been considered for solving
time-dependent linear kinetic transport equations, which arise for example
in electron radiation therapy or radiative heat transfer problems.
Many of the most popular methods are \emph{moment methods}, also known as \emph{moment closures} because they are distinguished by how they close the truncated
system of exact moment equations.
Moments are defined by angular averages against basis functions to
produce spectral approximations in the angle variable.
A typical family of moment models are the so-called \emph{$\PN$ methods}
\cite{Lewis-Miller-1984,Gel61}, which are pure spectral methods.
However, many high-order moment methods, including $\PN$, do not take
into account that the original kinetic particle distribution to be approximated must be
non-negative.
The moment vectors produced by such models are therefore often not realizable,
that is the fact that there is no associated non-negative kinetic distribution consistent with
the moment vector. Thus, the solutions can contain non-physical
artefacts such as negative local particle densities \cite{Bru02}.

The family of \emph{minimum-entropy models}, colloquially known as \emph{$\MN$ models}
or entropy-based moment closures, solve this problem (for certain physically
relevant entropies) by specifying the closure using a non-negative density
reconstructed from the moments. Additionally, they are hyperbolic and
dissipate entropy \cite{Levermore1998}.\\

All these properties of the minimum-entropy ansatz (for all types of angular bases) come at the price that the reconstruction of the distribution function involves solving an optimization problem at every point on the space-time mesh.
These reconstructions can be parallelized, and so the recent emphasis
on algorithms that can take advantage of massively parallel computing
environments has led to renewed interest in the computation of $\MN$ solutions
both for linear and non-linear kinetic equations
\cite{DubFeu99,Hauck2010,Lam2014,AllHau12,Garrett2014,
McDonald2012}.

To avoid these expensive computations a new class of models has been introduced in \cite{Ker76}, which will be called Kershaw closures. These models aim at providing a closed flux function which is generated by a non-negative distribution, thus ensuring that this crucial property of the transport solution is not violated. One derivation of Kershaw closures has been provided in \cite{Monreal}. Unfortunately, it is not straight forward to lift this procedure to higher dimensions. Here, a new ansatz is derived for which it is possible to do so.

This paper is organized as follows. First, the transport equation and its moment approximations are given. Then, the available realizability theory is shortly reviewed. Using these information allows to derive and investigate the class of Kershaw closures which is then intensively tested in well-known benchmark problems. Finally, conclusions and an outlook on future work is given.

\section{Modelling}
In slab geometry, the transport equation under consideration has the form 
\begin{align}
\label{eq:TransportEquation1D}
\dt\distribution+\SCheight\dz\distribution + \absorption\distribution = \scattering\collision{\distribution}+\source, \qquad \timevar\in\timeint,\z\in\Domain,\SCheight\in[-1,1].
\end{align}
The physical parameters are the absorption and scattering coefficient $\absorption,\scattering:\timeint\times\Domain\to\Rpos$, respectively, and the emitting source $\source:\timeint\times\Domain\times[-1,1]\to\Rpos$. Furthermore, $\SCheight\in[-1,1]$, and $\distribution = \distribution(\timevar,\z,\SCheight)$. 

The shorthand notation $\ints{\cdot} = \int\limits_{-1}^1\cdot~d\SCheight$ denotes integration over $[-1,1]$.

\begin{assumption}
\label{ass:CollisionOperator}
Following \cite{Levermore1996}, the collision operator $\collisionop$ is assumed to have the following properties.
\begin{enumerate}
\begin{subequations}
\label{eq:CollisionProperty}
\item Mass conservation
\begin{align}
\label{eq:CollisionPropertyMass}
\ints{\collision{\distribution}}=0.
\end{align}
\item Local entropy dissipation
\begin{align}
\label{eq:CollisionPropertyLocalDissipation}
\ints{\entropy'(\distribution)\collision{\distribution}}\leq 0,
\end{align}
where $\entropy$ denotes a strictly convex, twice differentiable entropy.
\item Constants in the kernel: 
\begin{align}
\label{eq:ConstantKernel}
\collision{c} = 0 \qquad \text{for every } c\in\R.
\end{align}
\end{subequations}
\end{enumerate}
\end{assumption}

{A typical example for $\collisionop$ is the linear integral operator
\begin{equation}
 \lincollision{\distribution} =  \int\limits_{-1}^1 \collisionkernel(\SCheight, \SCheight^\prime)
  \distribution(\timevar, \z, \SCheight^\prime)~d\SCheight^\prime 
  - \distribution(\timevar, \z, \SCheight),
\label{eq:collisionOperatorLin1D}
\end{equation}
where $\collisionkernel$ is non-negative, symmetric in both arguments and normalized to $\int\limits_{-1}^1 \collisionkernel(\SCheight, \SCheight^\prime)~d\SCheight^\prime=1$.
It includes the often-used special case of the BGK-type isotropic-scattering operator if $\collisionkernel\equiv \frac12$, which will be used here as well.}

\eqref{eq:TransportEquation1D} is supplemented by initial and boundary conditions:
\begin{subequations}
\begin{align}
\distribution(0,\z,\SCheight) &= \distributiontzero(\z,\SCheight) &\text{for } \z\in\Domain = (\zL,\zR), \SCheight\in[-1,1], \label{eq:TransportEquation1DIC}\\
\distribution(\timevar,\zL,\SCheight) &= \distributionboundary(\timevar,\zL,\SCheight) &\text{for } \timevar\in\timeint, \SCheight>0,  \label{eq:TransportEquation1DBCa}\\
\distribution(\timevar,\zR,\SCheight) &= \distributionboundary(\timevar,\zR,\SCheight) &\text{for } \timevar\in\timeint, \SCheight<0. \label{eq:TransportEquation1DBCb}
\end{align}
\end{subequations}

%Integrals over the ``half-sphere'' intervals $[-1,0]$ and $[0,1]$ are denoted by 
%\begin{align*}
%\intm{\cdot} = \int\limits_{-1}^{0}\cdot~d\SCheight,
%&&\intp{\cdot} = \int\limits_{0}^{1}\cdot~d\SCheight.
%\end{align*}
\section{Moment models and realizability}
In general, solving equation \eqref{eq:TransportEquation1D} is very expensive in higher dimensions due to the high dimensionality of the state space. 

For this reason it is convenient to use some type of spectral or Galerkin method to transform the high-dimensional equation into a system of lower-dimensional equations. Typically, one chooses to reduce the dimensionality by representing the angular dependence of $\distribution$ in terms of some basis $\basis$.
\begin{definition}
The vector of functions $\basis = \basis[\momentorder]:[-1,1]\to\R^{\momentorder+1}$ consisting of $\momentorder+1$ basis functions $\basiscomp[\basisind]$, $\basisind=0,\ldots\momentorder$ of maximal \emph{order} $\momentorder$ is called an \emph{angular basis}.

The so-called \emph{moments} of a given distribution function $\distribution$ with respect to $\basis$ are then defined by
\begin{align}
\label{eq:moments}
{\moments[] =\ints{{\basis}\distribution} = \left(\momentcomp{0},\ldots,\momentcomp{\momentorder}\right)^T},
\end{align}
where the integration is performed componentwise.\\

Assuming for simplicity $\basiscomp[0]\equiv 1$, the quantity $\momentcomp{0} = \ints{\basiscomp[0]\distribution}=\ints{\distribution}$ is called \emph{local particle density}. 
Furthermore, \emph{normalized moments} $\normalizedmoments = \left(\normalizedmomentcomp{1},\ldots,\normalizedmomentcomp{\momentorder}\right)\in\R^{\momentorder}$ are defined as 
\begin{align}
\label{eq:NormalizedMoments}
\normalizedmoments = \cfrac{\ints{\normalizedbasis\distribution}}{\ints{\distribution}}~,
\end{align}
where $\normalizedbasis = \left(\basiscomp[1],\ldots,\basiscomp[\momentorder]\right)^T$ is the remainder of the basis $\basis$.
\end{definition}
To obtain a set of equations for $\moments$, \eqref{eq:TransportEquation1D} has to be multiplied through by $\basis$ and integrated over $[-1,1]$, giving
\begin{align*}
\ints{\basis\dt\distribution}+\ints{\basis\dz\SCheight\distribution} + \ints{\basis\absorption\distribution} = \scattering\ints{\basis\collision{\distribution}}+\ints{\basis\source}.
\end{align*}
Collecting known terms, and interchanging integrals and differentiation where possible, the moment system has the form
\begin{align}
\label{eq:MomentSystemUnclosed1D}
\dt\moments+\dz\ints{\SCheight \basis\ansatz[\moments]} + \absorption\moments = \scattering\ints{\basis\collision{\ansatz[\moments]}}+\ints{\basis\source}.
\end{align}
The solution of \eqref{eq:MomentSystemUnclosed1D} is equivalent to the one of \eqref{eq:TransportEquation1D} if $\basis=\basis[\infty]$ is a basis of $\Lp{2}(\sphere,\R)$. 

Since it is impractical to work with an infinite-dimensional system, only a finite number of $\momentorder+1<\infty$ basis functions $\basis[\momentorder]$ of order $\momentorder$ can be considered. Unfortunately, there always exists an index $\basisind\in\{0,\dots,\momentorder\}$ such that the components of $\basiscomp\cdot\SCheight$ are not in the linear span of $\basis[\momentorder]$. Therefore, the flux term cannot be expressed in terms of $\moments[\momentorder]$ without additional information. Furthermore, the same might be true for the projection of the scattering operator onto the moment-space given by $\ints{\basis\collision{\distribution}}$. This is the so-called \emph{closure problem}. One usually prescribes some \emph{ansatz} distribution $\ansatz[\moments](\timevar,\spatialVariable,\SC):=\ansatz(\moments(\timevar,\spatialVariable),\basis(\SC))$ to calculate the unknown quantities in \eqref{eq:MomentSystemUnclosed1D}. Note that the dependence on the angular basis in the short-hand notation $\ansatz[\moments]$ is neglected for notational simplicity.\\

In this paper, the \emph{full-moment monomial basis} $\basiscomp = \SCheight^\basisind$ is considered. However, it is in principle possible to extend the derived concepts to other bases like half \cite{DubKla02,DubFraKlaTho03} or mixed moments \cite{Frank07,Schneider2014}.

\subsection{Minimum-entropy approach}
\label{sec:MinimumEntropy}\index{MinimumEntropy@\textbf{Minimum-entropy models}|(}
 An important class of models are the so-called \emph{minimum-entropy models}. Here the ansatz density $\ansatz$ is reconstructed from the moments $\moments$ by minimizing the \emph{entropy functional}
 \begin{align}
 \label{eq:entropyFunctional}
 \entropyFunctional(\distribution) = \ints{\entropy(\distribution)}
 \end{align}
 under the moment constraints \index{Entropy@\textbf{Entropy}!Entropy density $\entropy$}\index{Entropy@\textbf{Entropy}!Entropy functional $\entropyFunctional$}
 \begin{align}
 \label{eq:MomentConstraints}
 \ints{\basis\distribution} = \moments.
 \end{align}
The kinetic \emph{entropy density} $\entropy:\R\to\R$ is strictly convex and twice continuously differentiable and the minimum is simply taken over all functions $\distribution = \distribution(\SC)$ such that 
  $\entropyFunctional(\distribution)$ is well defined. The obtained ansatz $\ansatz = \ansatz[\moments]$ by solving the constrained optimization problem is given by\index{Ansatz@\textbf{Ansatz} $\ansatz$}
 \begin{equation}
  \ansatz[\moments] = \argmin\limits_{\distribution:\entropy(\distribution)\in\Lp{1}}\left\{\ints{\entropy(\distribution)}
  : \ints{\basis \distribution} = \moments \right\}.
 \label{eq:primal}
 \end{equation}
This problem is typically solved through its strictly convex finite-dimensional dual\index{Legendredual@\textbf{Legendre dual}}
 \begin{equation}
  \multipliers(\moments) := \argmin_{\tilde{\multipliers} \in \R^{\momentnumber}} \ints{\ld{\entropy}(\basis^T 
   \tilde{\multipliers})} - \moments^T \tilde{\multipliers},
 \label{eq:dual}
 \end{equation}
where $\ld{\entropy}$ is the Legendre dual \cite{courant2008methods} of $\entropy$.

Differentiating \eqref{eq:dual} with respect to the multipliers $\multipliers(\moments)$, the first-order necessary conditions show that the solution to \eqref{eq:primal} has the form
 \begin{equation}
  \ansatz[\moments] = \ld{\entropy}' \left(\basis^T \multipliers(\moments) \right),
 \label{eq:psiME}
 \end{equation}
where $\ld{\entropy}'$ is the derivative of $\ld{\entropy}$.\\

After substituting $\distribution$ in \eqref{eq:MomentSystemUnclosed1D} with $\ansatz[\moments]$, a closed system of equations remains, yielding
\begin{align}
\label{eq:MomentSystemClosed1D}
\dt\moments+\dz\cdot\ints{\SCheight\basis\ansatz[\moments]} + \absorption\moments = \scattering\ints{\basis\collision{\ansatz[\moments]}}+\ints{\basis\source}.
\end{align}
For convenience, \eqref{eq:MomentSystemClosed1D} can be written in the form of a usual first-order system of balance laws
\begin{align}
\label{eq:GeneralHyperbolicSystem1D}
\dt\moments+\dz\Flux_3\left(\moments\right) = \Source\left(\moments\right),
\end{align}
where 
\begin{subequations}
\label{eq:FluxDefinitions}
\begin{align}
\Flux\left(\moments\right) &= \ints{\SCheight\basis\ansatz[\moments]}\in\R^{\momentorder+1},\\
\Source\left(\moments\right) &= \scattering\ints{\basis\collision{\ansatz[\moments]}}+\ints{\basis\source}-\absorption\moments.
\end{align}
\end{subequations}
This system has several attractive properties, namely that it is hyperbolic and dissipates entropy \cite{Levermore1996,Levermore1998} and the eigenvalues of the flux Jacobian $\Flux'(\moments)$ are bounded in absolute value by one \cite{Schneider2015a}, which corresponds to the original velocities of the transport equation since $\SCheight\in[-1,1]$.

The kinetic entropy density $\entropy$ can be chosen according to the 
physics being modelled.
As in \cite{Levermore1996,Hauck2010}, \emph{Maxwell-Boltzmann entropy}
 \begin{align}
 \label{eq:EntropyM}
  \entropy(\distribution) = \distribution \log(\distribution) - \distribution
 \end{align}
is used, thus $\ld{\entropy}(p) = \ld{\entropy}'(p) = \exp(p)$. This entropy is used for non-interacting particles as in an ideal gas. Other physically relevant entropies are \cite{Levermore1996} 
\begin{align*}
  \entropy(\distribution) = \distribution \log(\distribution) + \left(1-\distribution\right)\log\left(1-\distribution\right)
\end{align*}
for particles satisfying \emph{Fermi-Dirac} (e.g. fermions) or
\begin{align*}
  \entropy(\distribution) = \distribution \log(\distribution) - \left(1+\distribution\right)\log\left(1+\distribution\right)
\end{align*}
for particles satisfying \emph{Bose-Einstein statistics} (e.g. bosons).

The resulting minimum-entropy model is commonly referred to as the $\MN$ model.

Note that the classical $\PN$ approximation \cite{Eddington,Lewis-Miller-1984} is also an 
entropy-based moment closure using the entropy
$$
  \entropy(\distribution) = \frac12 \distribution^2 = \ld{\entropy}(\distribution)
$$
and Legendre polynomials as angular basis \cite{Brunner2005b,Seibold2014}.

\subsection{Realizability}
Since the underlying kinetic density to be approximated is
non-negative, a 
moment vector only makes sense physically if it can be associated with a 
non-negative distribution function. In this case the moment vector is called 
\emph{realizable}.

\begin{definition}
\label{def:RealizableSet}
The \emph{realizable set} $\RD{\basis}{}$\index{Realizability@\textbf{Realizability}!Realizable set $\RD{\basis}{}$} is 
$$
\RD{\basis}{} = \left\{\moments~:~\exists \distribution(\SCheight)\ge 0,\, \ints{\distribution} > 0,
 \text{ such that } \moments =\ints{\basis\distribution} \right\}.
$$
If $\moments\in\RD{\basis}{}$, then $\moments$ is called \emph{realizable}.
Any $\distribution$ such that $\moments =\ints{\basis \distribution}$ is called a \emph{representing 
density}. If $\distribution$ is additionally a linear combination of Dirac deltas \cite{Hassani2009,Mathematics2011,Kuo2006}, it is called \emph{atomic} \cite{Curto1991}.
\end{definition}

\begin{definition}
\mbox{ }
\begin{enumerate}[(a)]
\item Let $\hankelA,\hankelB\in\R^{\momentnumber\times \momentnumber}$ be Hermitian matrices. The partial ordering $"\geq"$ on such matrices is defined by $\hankelA\geq \hankelB$ if and only if $\hankelA-\hankelB$ is positive semi-definite. In particular $\hankelA\geq 0$ denotes that $\hankelA$ is positive semi-definite.

\item A distribution function $\distribution$ is said to be $\hankelrank-$\emph{atomic} with \emph{atoms} $\SCheight_\hankelind$ and \emph{densities} $\hankelrho_\hankelind$ if it is a linear-combination of $\hankelrank$ Dirac deltas of the form
\begin{align*}
\distribution = \sum_{\hankelind=0}^{\hankelrank-1}\hankelrho_\hankelind\dirac(\SCheight-\SCheight_\hankelind).
\end{align*}
\end{enumerate}
\end{definition}

For the full-moment basis the question of finding practical characterizations of the realizable set $\RD{\basis}{}$ has been completely solved in \cite{Curto1991}. The resulting problems are special cases of the so-called \emph{truncated Hausdorff moment problem} which aims to find a realizing distribution with support on the interval $[a,b]$ for a moment vector $\moments$ in the set
$$\left\{\moments\in\R^{\momentorder+1}~|~ \exists \distribution\geq 0 \text{ with } \momentcomp{\basisind} = \int\limits_a^b \SCheight^\basisind\distribution~d\SCheight,~\basisind=0,\ldots,\momentorder\right\}.$$

The following characterizations of the above realizable set hold.
{\begin{lemma}[Truncated Hausdorff moment problem \cite{Curto1991}]
\label{lem:FullMomentRealizability}
%For the truncated Hausdorff moment problem, the following realizability conditions hold.
Define the \emph{Hankel matrices}
$$
\hankelA(\hankelhalfind):=\left(\momentcomp{i+j}\right)_{i,j=0}^\hankelhalfind, \quad \hankelB(\hankelhalfind):=\left(\momentcomp{i+j+1}\right)_{i,j=0}^\hankelhalfind, \quad \hankelC(\hankelhalfind):=\left(\momentcomp{i+j}\right)_{i,j=1}^\hankelhalfind.
$$
Then the truncated Hausdorff moment problem has a solution if and only if
\begin{itemize}
\item for $\momentorder=2\hankelhalfind+1$,
\begin{align}
%A(\hankelhalfind)\geq 0 \qquad \text{and} \\
b\hankelA(\hankelhalfind)\geq \hankelB(\hankelhalfind)\geq a\hankelA(\hankelhalfind); \label{eq:haus-odd}
\end{align}
\item for $\momentorder=2\hankelhalfind$,
\begin{align}
\hankelA(\hankelhalfind) &\geq 0 \qquad \text{and} \label{eq:haus-even-1} \\
(a+b)\hankelB(\hankelhalfind-1) &\geq ab\hankelA(\hankelhalfind-1) + \hankelC(\hankelhalfind). \label{eq:haus-even-2}
\end{align} 
\end{itemize}
\end{lemma}
}
\begin{corollary}
\label{cor:FullMomentRealizability}
%Define the \emph{Hankel matrices}
%$$
%\hankelA(\hankelhalfind):=\left(\momentcomp{i+j}\right)_{i,j=0}^\hankelhalfind, \quad \hankelB(\hankelhalfind):=\left(\momentcomp{i+j+1}\right)_{i,j=0}^\hankelhalfind, \quad \hankelC(\hankelhalfind):=\left(\momentcomp{i+j}\right)_{i,j=1}^\hankelhalfind.
%$$
The realizable set $\RD{\fmbasis}{}$ satisfies
\begin{align*}
\RD{\fmbasis}{} = \begin{cases}
\left\{\moments\in\R^{\momentorder+1}~|~ {\hankelA(\hankelhalfind)\geq \hankelB(\hankelhalfind),~\hankelA(\hankelhalfind)\geq -\hankelB(\hankelhalfind)}\right\} & \text{ if } $\momentorder=2\hankelhalfind+1$,\\
\left\{\moments\in\R^{\momentorder+1}~|~ \hankelA(\hankelhalfind)\geq 0, \hankelA(\hankelhalfind-1)\geq \hankelC(\hankelhalfind)\right\} & \text{ if } $\momentorder=2\hankelhalfind$.\\
\end{cases}
\end{align*}
\end{corollary}
\begin{proof}
{Follows directly from \lemref{lem:FullMomentRealizability} with $a = -1$, $b=1$.}
\end{proof}

Furthermore, it can be shown that there exists a minimal atomic representing distribution $\distribution$ (in the sense that it contains the fewest possible number of atoms while still representing the moments) and that one can directly find this distribution with the help of its \emph{generating function}. This is the consequence of a recursiveness property of the Hankel matrices $\hankelA(\hankelhalfind)$ and $\hankelB(\hankelhalfind)$ \cite{Curto1991}.

\begin{corollary}
\label{cor:GeneratingFunction}
Let $\left(\hankelphi_0,\ldots,\hankelphi_{\hankelrank-1}\right)^T:= \hankelA(\hankelrank-1)^{-1}\hankelv(\hankelrank,\hankelrank-1)$, where $\hankelv(i,j) = \left(\momentcomp{i+l}\right)_{l=0}^j$ and $\hankelrank$ is the smallest integer such that $\hankelA(\hankelrank)$ is singular,%
\footnote{$\hankelrank$ is also called the \emph{Hankel rank} \cite{Curto1991}.} %
the generating function is defined by
\begin{align*}
\generatingFunction(\SCheight) = \SCheight^\hankelrank-\sum\limits_{\hankelind=0}^{\hankelrank-1}\hankelphi_\hankelind\SCheight^\hankelind.
\end{align*}
The roots of this polynomial give the atoms $\SCheight_\hankelind$ of the distribution $\distribution = \sum_{\hankelind=0}^{\hankelrank-1}\hankelrho_\hankelind\dirac(\SCheight-\SCheight_\hankelind)$ whereas the densities $\hankelrho_\hankelind$ are calculated afterwards from the Vandermonde system 
\begin{align}
\label{eq:VandermondeSystem}
\rho_0\SCheight_0^\hankelind+\cdots +\rho_{\hankelrank-1}\SCheight_{\hankelrank-1}^\hankelind = \hankelv(\hankelind,0) = \momentcomp{\hankelind}~~~\hankelind=0\ldots \hankelrank-1.
\end{align}
Furthermore, when the moment vector $\moments$ is on the boundary of the realizable set, the minimal atomic representing measure is the unique representing measure.
\end{corollary}
%\begin{remark}
%Note that finding the atoms $\SCheight_i$ can be done in a robust way. Although the condition number of the Vandermonde system is exponential in the Hankel rank $\hankelrank$ \cite{Gautschi1987}, efficient and robust algorithms for the solution of this linear system exist \cite{Bjorck1970}. However, all results presented here are obtained analytically using symbolic calculations.
%\end{remark}
Due to the structure of the used Hankel matrices (the highest moment $\momentcomp{\momentorder}$ always appears exactly once in the entries of the matrices) it is always possible to rearrange the conditions involving this highest moment in \corref{cor:FullMomentRealizability} in such a way that 
\begin{align*}
%\label{eq:UpperLowerBoundsFullMoments}
\fup(\momentcomp{0},\ldots,\momentcomp{\momentorder-1})\geq \momentcomp{\momentorder} \geq \flow(\momentcomp{0},\ldots,\momentcomp{\momentorder-1})
\end{align*}
for functions $\fup$ and $\flow$. Whenever $\moments$ is realizable and $\momentcomp{\momentorder} = \flow(\momentcomp{0},\ldots,\momentcomp{\momentorder-1})$, $\moments$ is said to be on the \emph{lower $\momentorder^{\text{th}}$-order realizability boundary}. Similarly, if $\momentcomp{\momentorder} = \fup(\momentcomp{0},\ldots,\momentcomp{\momentorder-1})$, $\moments$ is said to be on the \emph{upper $\momentorder^{\text{th}}$-order realizability boundary}.

The functions $\fup$ and $\flow$ can be specified using the pseudoinverses of $\hankelA(\hankelhalfind-1)$ and $\hankelA(\hankelhalfind-1)\pm \hankelB(\hankelhalfind-1)$ due to Lemma 2.3 in \cite{Curto1991}. Also compare \cite{Schneider2014} where this relation is used for the mixed-moment setting.

{
\begin{lemma}[cf. Lemma 2.3 in \cite{Curto1991}]
\label{lem:pd-ext}
\mbox{ }\\
Let $\hankelA \in \R^{\hankelhalfind \times\hankelhalfind}$ be symmetric, $\hankelb \in \R^\hankelhalfind$, and $c \in \R$ such that
$$
\tilde\hankelA = \left(\begin{array}{cc} \hankelA & \hankelb \\ \hankelb^T & c \end{array}\right).
$$
\begin{itemize}
 \item[(i)] If $\tilde\hankelA \ge 0$, then $\hankelA \ge 0$, $\hankelb = \hankelA \hankelw \in \range{\hankelA}$ and $c \ge \hankelw^T \hankelA \hankelw$.
 \item[(ii)] If $\hankelA \ge 0$ and $\hankelb = \hankelA \hankelw\in \range{\hankelA}$, then $\tilde\hankelA \ge 0$ if and only if $c \ge \hankelw^T \hankelA \hankelw$.
\end{itemize}
\end{lemma}
}
\begin{remark}
\label{rem:Pseudoinverse}
Since $\hankelA$ is invertible on its range, a more explicit formula for the bound on $c$ in Lemma \ref{lem:pd-ext} can be written using the pseudo-inverse $\pseudoinv{\hankelA}$ of $\hankelA$. That is, for all $\hankelw$ such that $\hankelb = \hankelA \hankelw$, it is also $\hankelw^T \hankelA \hankelw = \hankelb^T \pseudoinv{\hankelA} \hankelb$.  Thus the bound on $c$ is well-defined even when $\hankelA$ is singular.
\end{remark}

To simplify notation later, the following corollary is written in terms of $\momentcomp{\momentorder+1}$ instead of $\momentcomp{\momentorder}$.
\begin{corollary}
\label{cor:FupFlow}
The functions $\fup$ and $\flow$ satisfying 
\begin{align}
\label{eq:UpperLowerBoundsFullMoments}
\fup(\momentcomp{0},\ldots,\momentcomp{\momentorder})\geq \momentcomp{\momentorder+1} \geq \flow(\momentcomp{0},\ldots,\momentcomp{\momentorder})
\end{align}
are given by
\begin{align*}
\fup(\momentcomp{0},\ldots,\momentcomp{\momentorder}) &= 
\begin{cases}
\momentcomp{\momentorder-1}-\hankelb_-^T\pseudoinv{\left(\hankelA(\hankelhalfind-1)-\hankelC(\hankelhalfind-1)\right)}\hankelb_- & \text{ if } \momentorder = 2\hankelhalfind+1\\
\momentcomp{\momentorder}-\hankelb_-^T\pseudoinv{\left(\hankelA(\hankelhalfind-1)-\hankelB(\hankelhalfind-1)\right)}\hankelb_- & \text{ if } \momentorder = 2\hankelhalfind
\end{cases}\\
\flow(\momentcomp{0},\ldots,\momentcomp{\momentorder}) &= 
\begin{cases}
\hankelb_+^T\pseudoinv{\hankelA}(\hankelhalfind)\hankelb_+& \text{ if } \momentorder = 2\hankelhalfind+1\\
-\momentcomp{\momentorder}+\hankelb_+^T\pseudoinv{\left(\hankelA(\hankelhalfind-1)+\hankelB(\hankelhalfind-1)\right)}\hankelb_+~ & \text{ if } \momentorder = 2\hankelhalfind
\end{cases}
\end{align*}
where in the odd case
\begin{align*}
\hankelb_- = \left(\momentcomp{\hankelhalfind}-\momentcomp{\hankelhalfind+2},\ldots,\momentcomp{\momentorder-2}-\momentcomp{\momentorder}\right)^T,\qquad \hankelb_+ = \left(\momentcomp{\hankelhalfind+1},\ldots,\momentcomp{\momentorder}\right)^T
\end{align*}
and in the even case
\begin{align*}
\hankelb_\mp = \left(\momentcomp{\hankelhalfind}\mp\momentcomp{\hankelhalfind+1},\ldots,\momentcomp{\momentorder-1}\mp\momentcomp{\momentorder}\right)^T.
\end{align*}
\end{corollary}
\begin{proof}
For convenience, a proof of the even case if given. The odd case follows similarly.

If $\momentorder = 2\hankelhalfind>0$, the $(\momentorder+1)^\text{th}$-order realizable set is given by 
\begin{align*}
\RD{\fmbasis+1}{} =
\left\{\moments[\momentorder+1]\in\R^{\momentorder+2}~|~ {\hankelA(\hankelhalfind)\geq \hankelB(\hankelhalfind),~\hankelA(\hankelhalfind)\geq -\hankelB(\hankelhalfind)}\right\}.
\end{align*}
The moment constraints have the form 
\begin{align*}
\hankelA(\hankelhalfind)\mp\hankelB(\hankelhalfind) = 
\begin{pmatrix}
\hankelA(\hankelhalfind-1)\mp\hankelB(\hankelhalfind-1) & \hankelb_\mp\\\hankelb_\mp^T & \momentcomp{\momentorder}\mp \momentcomp{\momentorder+1}
\end{pmatrix}\geq 0.
\end{align*}
Thus, by \lemref{lem:pd-ext} and \rmref{rem:Pseudoinverse}, it holds that
\begin{align}
\label{eq:bounds1}
\momentcomp{\momentorder}\mp \momentcomp{\momentorder+1} &\geq \hankelb_\mp^T\pseudoinv{\left(\hankelA(\hankelhalfind-1)\mp\hankelB(\hankelhalfind-1)\right)}\hankelb_\mp = \hankelw_\mp^T\left(\hankelA(\hankelhalfind-1)\mp\hankelB(\hankelhalfind-1)\right)\hankelw_\mp~,\\
\label{eq:bounds2}
\left(\hankelA(\hankelhalfind-1)\mp\hankelB(\hankelhalfind-1)\right)\hankelw_\mp &= \hankelb_\mp
\end{align}
Consequently,
\begin{align*}
\momentcomp{\momentorder}-\hankelb_-^T\pseudoinv{\left(\hankelA(\hankelhalfind-1)-\hankelB(\hankelhalfind-1)\right)}\hankelb_- 
\geq \momentcomp{\momentorder+1} \geq
-\momentcomp{\momentorder}+\hankelb_+^T\pseudoinv{\left(\hankelA(\hankelhalfind-1)+\hankelB(\hankelhalfind-1)\right)}\hankelb_+~.
\end{align*}
\end{proof}
\begin{remark}
By convention, $\hankelb_-^T\pseudoinv{\left(\hankelA(\hankelhalfind-1)-\hankelC(\hankelhalfind-1)\right)}\hankelb_- = 0$ if $\momentorder = 1$.
\end{remark}
{
\begin{remark}
An equivalent description for $\fup$ and $\flow$ has been derived in \cite{Ker76} using determinants of the Hankel matrices. However, that description is only well-defined if these matrices are invertible. This can lead to numerical problems close to the realizability boundary.
In implementation, the formulation given in \lemref{lem:pd-ext} proved to be the most stable. For example in the proof of \corref{cor:FupFlow} it is more stable to first calculate the solution $\hankelw_\mp$ of \eqref{eq:bounds2} (this can be handled efficiently even in the singular case using the \emph{backslash} operator in Matlab \cite{MATLAB:2012}) and plug it into the quadratic form on the right-hand side of \eqref{eq:bounds1} instead of calculating the pseudoinverse of $\hankelA(\hankelhalfind-1)\mp\hankelB(\hankelhalfind-1)$.
\end{remark}
}

\begin{example}
\label{ex:M2Realizability}
For the full-moment setting $\fmbasis[2]$, i.e. $\hankelhalfind=1$, the Hankel matrices have the form
\begin{align*}
\hankelA(0) = \momentcomp{0},~
\hankelB(0) = \momentcomp{1},~
\hankelC(0) = \momentcomp{2},~
\hankelA(1) = \begin{pmatrix}
\momentcomp{0} & \momentcomp{1}\\
\momentcomp{1} & \momentcomp{2}
\end{pmatrix}.
\end{align*}
Sylvester's criterion implies that a symmetric matrix is positive semi-definite if its leading principal minors are non-negative \cite{meyer2000matrix}. Thus, all moments in $\RD{\basis[2]}{}$ satisfy
\begin{align}
\label{eq:M2Realizability}
\momentcomp{0}\geq 0, \qquad\momentcomp{0}\momentcomp{2}\geq \momentcomp{1}^2,\quand \momentcomp{0}\geq\momentcomp{2}.
\end{align}
Note that these conditions also imply the first-order realizability condition 
\begin{align*}
\pm\momentcomp{1}\leq \momentcomp{0}.
\end{align*}
The lower second-order realizability boundary is given by $\momentcomp{0}\momentcomp{2} = \momentcomp{1}^2$ while $\momentcomp{0}=\momentcomp{2}$ defines the upper one.\\

The third-order realizability conditions in the non-singular case $\pm \momentcomp{1} < \momentcomp{0}$ are given by \cite{Ker76}
\begin{align}
\label{eq:FullMomentsThirdOrderRealizabilityConditions}
\momentcomp{2}-\cfrac{(\momentcomp{1}-\momentcomp{2})^2}{\momentcomp{0}-\momentcomp{1}}\geq \momentcomp{3} \geq -\momentcomp{2}+\cfrac{(\momentcomp{1}+\momentcomp{2})^2}{\momentcomp{0}+\momentcomp{1}}~.
\end{align}
\end{example}

\begin{remark}
Due to the structure of the Hankel matrices, the upper and lower bounds can be expressed in terms of normalized moments $\normalizedmoments$ as
\begin{align}
\label{eq:UpperLowerBoundsFullMomentsNorm}
\fup(1,\ldots,\normalizedmomentcomp{\momentorder})\geq \normalizedmomentcomp{\momentorder+1} \geq \flow(1,\ldots,\normalizedmomentcomp{\momentorder}).
\end{align}
\end{remark}
\section{Kershaw closures}
\subsection{Derivation}
With the previous realizability theory it is now possible to develop another closure strategy, which is called \emph{Kershaw} closure. The key idea of Kershaw in \cite{Ker76} was to derive an easy closure that preserves the realizability conditions, i.e. for every $\moments[{\basis[\momentorder]}]\in \RD{\basis[\momentorder]}{}$ the moment vector including the higher-order moments of the flux, generated by the closure relation, satisfies $\moments[{\basis[\momentorder+1]}]\in \RD{\basis[\momentorder+1]}{}$. If one uses the non-negative representing (atomic) distributions provided by the realizability theorems above, this is satisfied automatically. 
A very important property is defined as follows.
\begin{definition}
\mbox{ }\\
A representing distribution $\distribution$ for $\moments$ is said to \emph{interpolate the isotropic point} if
$$\ints{\basis[\momentorder]\distribution} = \frac12\ints{\distribution}\ints{\basis[\momentorder]} = \frac12\ints{\distribution}\isotropicmoment ~\text{ implies }~ \ints{\basis[\momentorder+1]\distribution} = \frac12\ints{\distribution}\ints{\basis[\momentorder+1]},$$
i.e. the moment vector including the $(\momentorder+1)^{\text{th}}$ moment is also isotropic with respect to $\basis[\momentorder+1]$.
\end{definition}

Unfortunately, it is in general wrong that this atomic distribution correctly interpolates the isotropic point, but it is possible to cure this problem abusing the structure of the bounded realizable set. 

Since $\RDone{\basis[\momentorder+1]}{} = \left\{\moments\in \RD{\basis[\momentorder+1]}{}~|~\momentcomp{0} = 1 \right\}$ is convex and bounded \cite{Schneider2014}, every vector $\moments[{\basis[\momentorder+1]}]\in\RDone{\basis[\momentorder+1]}{}$ can be written as a convex combination of moment vectors on the realizability boundary. Since all these boundary moments are by definition realizable, and due to the linearity of the moment problem, the convex combination of their representing densities is a representing density for $\moments[{\basis[\momentorder+1]}]$. 

All that remains is to choose the convex combination in such a way that the isotropic point is correctly interpolated. However, this requires $(\momentorder+1)^{\text{th}}$-order realizability information for a $\momentorder^{\text{th}}$-order closure.

%\begin{remark}
%\mbox{ }
%\begin{enumerate}[(a)]
%\item Having this interpolation property is very crucial for a moment system. If it is not satisfied, the corresponding moment system cannot converge to the correct solution in the diffusive (or isotropic) limit (compare \secref{sec:ApplicationDiffusionApproximation} and \secref{sec:ApplicationMomentApproximation}). Even more, since the constants are in the kernel of the collision operators $\lincollisionop$ and $\LaplaceBeltramiProjection$, the moment system cannot converge to the correct equilibrium state (with respect to the transport equation), if this state is homogeneous in space.
%
%\item Choosing the closure in such a way also ensures that the solution of the moment system is exact (with respect to the transport part of the transport equation) if the moments of the true transport solution are on the parts of the boundary of the corresponding realizable set where the atomic representing distribution is unique (recall \corref{cor:GeneratingFunction} and \rmref{rem:MMnNonUniquenessOnRealizabilityBoundary}). This is sometimes called the \emph{anisotropic free-streaming limit}.
%
%\item Since the underlying representing distribution is atomic, the evaluation of point-wise values (in $\SCheight$) does not make sense. This is in particular a problem for mixed moments with the Laplace-Beltrami operator, where the distribution has to be evaluated at $\SCheight=0$ (compare \secref{sec:MixedMoments1D}).
%\end{enumerate}
%\end{remark}

\begin{figure}[h]
\externaltikz{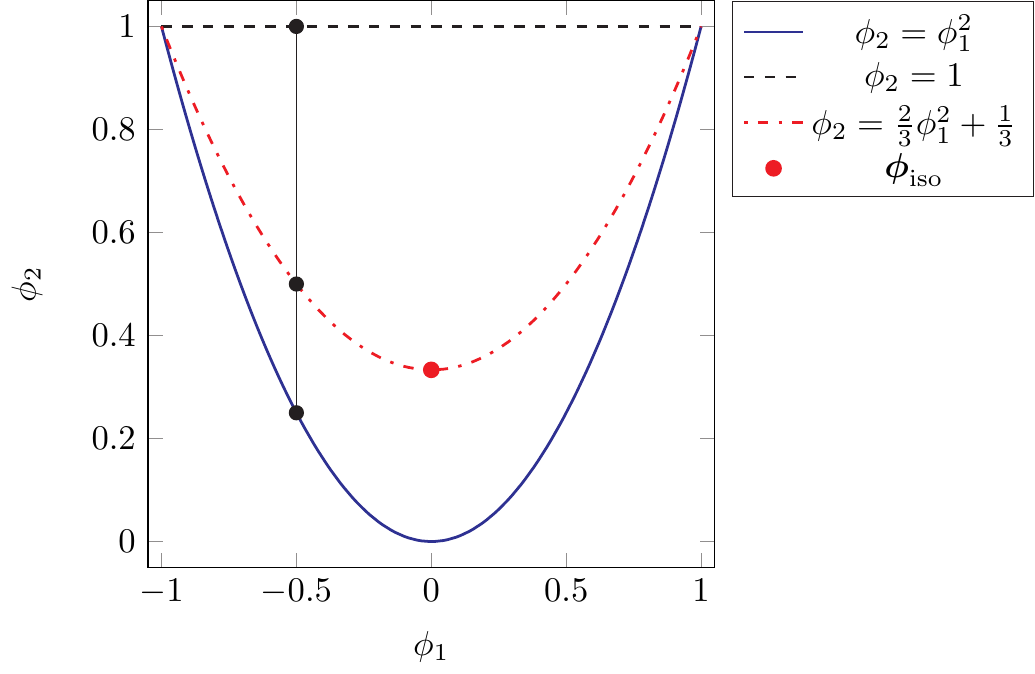}{\relinput{Images/K1Construction}}
 \centering
 \caption{Construction of the $\KN[1]$ closure. The red dash-dotted curve is the convex combination of the upper (black, dashed) and lower (blue, solid) second-order boundary-moments which interpolates the second-order isotropic point (red dot).}
 \label{fig:K1Construction}
\end{figure}

Since the above definition of a Kershaw closure is very abstract, the procedure is introduced in detail for $\fmbasis[1]$. Recall the second-order realizability conditions \eqref{eq:M2Realizability}. In normalized moments, they are equivalent to
$$
-1 \leq \normalizedmomentcomp{1}\leq 1,\quad \normalizedmomentcomp{1}^2\leq \normalizedmomentcomp{2}\leq 1.
$$
Thus the normalized second moment $\normalizedmomentcomp{2}$ is bounded from above and below by $\fup(\normalizedmomentcomp{1}) = 1$ and $\flow(\normalizedmomentcomp{1}) = \normalizedmomentcomp{1}^2$ depending only on $\normalizedmomentcomp{1}$.
The distribution on the lower second-order realizability boundary ($\normalizedmomentcomp{2} = \normalizedmomentcomp{1}^2$) is given by 
\begin{align}
\label{eq:K1psilow}
\ansatz[\text{low}] = \momentcomp{0}\dirac\left(\normalizedmomentcomp{1}-\SCheight\right)
\end{align}
while the upper second-order realizability-boundary distribution ($\normalizedmomentcomp{2}=1$) is given by
\begin{align}
\label{eq:K1psiup}
\ansatz[\text{up}] = \momentcomp{0}\left(\frac{1-\normalizedmomentcomp{1}}{2}\dirac(1+\SCheight) + \frac{1+\normalizedmomentcomp{1}}{2}\dirac(1-\SCheight)\right).
\end{align} 
By linearity of the problem, every convex combination 
\begin{align}
\label{eq:K1Distribution}
\ansatz = \convexscalar\ansatz[\text{low}] + (1-\convexscalar)\ansatz[\text{up}],\quad \convexscalar\in[0,1]
\end{align}
reproduces all moments up to first order and satisfies $\ansatz\geq 0$. 
\begin{remark}
The choice of $\convexscalar = \convexscalar(\normalizedmoments)$ is completely arbitrary. There are only two conditions, namely $\convexscalar\in[0,1]$ (convex combination property) and $\convexscalar(\normalizedisotropicmoment)$ is chosen such that the isotropic point is correctly interpolated. If, for whatever reason, other points should be interpolated as well, these interpolation conditions can be satisfied similarly.

For simplicity, $\convexscalar$ is chosen to be constant. This also ensures that the described procedure is uniquely determined.
\end{remark}
Calculating normalized isotropic moments up to order two gives $\normalizedisotropicmoment = (0,\frac{1}{3})$. The normalized second moment of \eqref{eq:K1Distribution} at the isotropic point therefore satisfies
\begin{align*}
\normalizedmomentcomp{2} = \convexscalar (\normalizedmomentcomp{1})^2+(1-\convexscalar) \stackrel{\normalizedmomentcomp{1}=0}{=} (1-\convexscalar) \stackrel{!}{=} \frac{1}{3}.
\end{align*}
This implies $\convexscalar = \frac{2}{3}$, altogether resulting in the analytically closed directive for the normalized second moment
\begin{align}
\label{eq:K1closure}
\normalizedmomentcomp{2} = \frac23\normalizedmomentcomp{1}^2+\frac13.
\end{align}
This construction procedure is shown in \figref{fig:K1Construction}, where the corresponding upper and lower realizability conditions and their convex combination interpolating the isotropic point are plotted.\\

This model is called the \emph{Kershaw $\KN[1]$ closure}. Similarly, higher-order Kershaw closures can be derived. If only the flux has to be closed, it suffices to calculate the interpolation on the moment level.

\begin{figure}[h]
\settikzlabel{fig:K2Flux}\settikzlabel{fig:K2FluxDifference}
\centering
\externaltikz{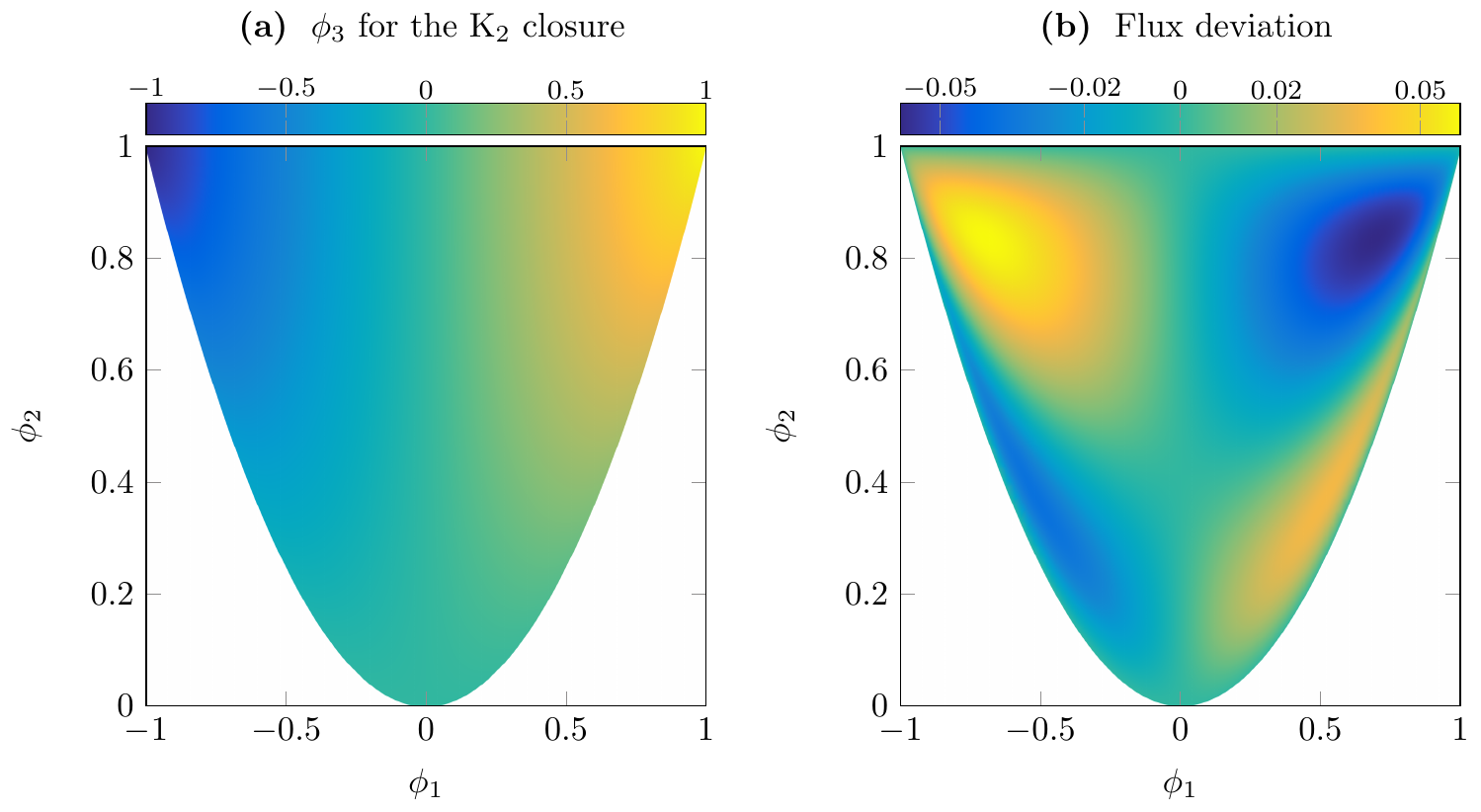}{\relinput{Images/K2FluxAndDifference}}
\caption{Third normalized moment $\normalizedmomentcomp{3}$, as in \eqref{eq:K2closure}, for the $\KN[2]$ model and its deviation from the corresponding $\MN[2]$ moment.}
\label{fig:K2FluxAndDifference}
\end{figure}

\begin{example}[${\KN[2]}$ closure]
The third-order realizability conditions \eqref{eq:FullMomentsThirdOrderRealizabilityConditions} can be written as
\begin{align*}
\normalizedmomentcomp{2}-\cfrac{(\normalizedmomentcomp{1}-\normalizedmomentcomp{2})^2}{1-\normalizedmomentcomp{1}}\geq \normalizedmomentcomp{3} \geq -\normalizedmomentcomp{2}+\cfrac{(\normalizedmomentcomp{1}+\normalizedmomentcomp{2})^2}{1+\normalizedmomentcomp{1}}~.
\end{align*}
The normalized isotropic moment of third order is $\normalizedisotropicmoment = \left(0,\frac13,0\right)^T$. To obtain the $\KN[2]$ closure, the ansatz is again a convex combination of upper and lower boundary moments, i.e.
\begin{align}
\label{eq:Kershaw2ConvexCombination}
\normalizedmomentcomp{3}(\normalizedmoments) = \convexscalar\left(-\normalizedmomentcomp{2}+\cfrac{(\normalizedmomentcomp{1}+\normalizedmomentcomp{2})^2}{1+\normalizedmomentcomp{1}}\right)
+(1-\convexscalar)\left(\normalizedmomentcomp{2}-\cfrac{(\normalizedmomentcomp{1}-\normalizedmomentcomp{2})^2}{1-\normalizedmomentcomp{1}}\right).
\end{align}
Evaluating \eqref{eq:Kershaw2ConvexCombination} at $\normalizedmomentcomp{1} = 0$, $\normalizedmomentcomp{2} = \frac13$ gives
\begin{align*}
\normalizedmomentcomp{3}(\normalizedisotropicmoment) = \frac{4\, \convexscalar}{9} - \frac{2}{9} \stackrel{!}{=} 0,
\end{align*}
which implies $\convexscalar = \frac12$. Therefore, the $\KN[2]$ closure for the third moment is given by
\begin{align}
\label{eq:K2closure}
\normalizedmomentcomp{3}(\normalizedmoments) = \frac{\normalizedmomentcomp{1}\, \left(\normalizedmomentcomp{1}^2 + \normalizedmomentcomp{2}^2 - 2\, \normalizedmomentcomp{2}\right)}{\normalizedmomentcomp{1}^2 - 1}.
\end{align}
See also \figref{fig:K2FluxAndDifference} for a visualization of \eqref{eq:K2closure} and its slight deviation from the minimum-entropy $\MN[2]$ model.
\end{example}

\begin{remark}
Note that the interpolation procedure given in \cite{Monreal} produces different results for $\momentorder>1$. There, the $\KN[2]$ closure is given by 
\begin{align*}
\normalizedmomentcomp{3}(\normalizedmoments) = \normalizedmomentcomp{1}\normalizedmomentcomp{2}.
\end{align*}
{
These two representations are identical on the realizability boundary but differ slightly in the interior. However, the overall structure generated by this approach is completely different, as investigated later on in \secref{sec:Eigenstructure}. Furthermore, the here-derived model is consistent with the construction of the mixed-moment Kershaw closures in \cite{Schneider2014}.
}
\end{remark}

{
\begin{theorem}
\mbox{ }\\
The \emph{Kershaw closure} $\KN$ of order $\momentorder$ is given by
\begin{align}
\label{eq:KnAnsatz}
\normalizedmomentcomp{\momentorder+1}(\normalizedmoments) = \convexscalar \flow(\normalizedmoments)+(1-\convexscalar) \fup(\normalizedmoments),
\end{align}
where the interpolation constant
\begin{align}
\label{eq:KnAnsatzScalar}
\convexscalar = \cfrac{\frac12\ints{\SCheight^{\momentorder+1}}-\fup(\normalizedisotropicmoment)}{\flow(\normalizedisotropicmoment)-\fup(\normalizedisotropicmoment)} = 
\begin{cases}
\frac{\hankelhalfind+2}{2\hankelhalfind+3} & \text{ if } \momentorder = 2\hankelhalfind+1\\
\frac12 & \text{ if } \momentorder = 2\hankelhalfind
\end{cases}
\end{align}
is defined via the functions $\fup$ and $\flow$ as given in \corref{cor:FupFlow}.
\end{theorem}
}
\begin{proof}
Using the ansatz \eqref{eq:KnAnsatz}, the interpolation condition yields
\begin{align*}
\normalizedmomentcomp{\momentorder+1}(\normalizedisotropicmoment) = \convexscalar \flow(\normalizedisotropicmoment)+(1-\convexscalar) \fup(\normalizedisotropicmoment) \stackrel{!}{=} \frac12\ints{\SCheight^{\momentorder+1}}.
\end{align*}
Solving this equation for $\convexscalar$ gives the first equality in \eqref{eq:KnAnsatzScalar}. Since the proof of the second equality only requires simple algebra and analysis but is tedious and unenlightening, it is omitted here.
\end{proof}

\begin{remark}
The big advantage of Kershaw closures is that the calculation of the fluxes is very cheap compared to minimum-entropy models (recall that they require to solve the non-linear moment system \eqref{eq:dual} to calculate the flux) while mimicking their behaviour. Indeed, the more robust $\hankelw$-notation in \corref{cor:FupFlow} only requires the solution of four linear systems of dimension at most $\hankelhalfind+1$. Even more, for $\momentnumber$ not too big, the closures can be derived symbolically leading to an even greater gain in efficiency.

{Also note that the class of Kershaw closures is related to the often-used \emph{quadrature method of moments}. See \cite{Schneider2014} for a discussion of the similarities and differences between these two classes of models.}
\end{remark}

\subsection{Eigenvalues and characteristic fields}
\label{sec:Eigenstructure}

While much is known about the eigenstructure of minimum-entropy models, until now no general theory is available for Kershaw closures. This section briefly investigates the eigenstructure and the characteristic fields of the first two models in this hierarchy to show that they are indeed hyperbolic.

Due to the structure of the flux, similar to the minimum-entropy equations, the Jacobian $\Flux'$ always has ones on the first super-diagonal and last row $\left(\nabla_{\moments[]} \momentcomp{\momentorder+1}\right)^T$, such that every eigenvalue $\lambda$ of $\Flux'$ has the eigenvector
\begin{align}
\label{eq:MnEigenvectors}
\eigenvector &= \left(1,\eigenvalue,\eigenvalue^2,\ldots,\eigenvalue^{\momentorder}\right)^T\quad\text{ with }\\
\label{eq:MnEigenvalueEq}
0 &= \eigenvalue^{\momentorder+1}-\left(\nabla_{\moments[]} \momentcomp{\momentorder+1}\right)^T\eigenvector.
\end{align}

\begin{example}[The {$\KN[1]$} closure]
Given \eqref{eq:K1closure}, it is easy to conclude that
\begin{align}
\Flux'(\moments) = \left(\begin{array}{cc} 0 & 1\\ \frac{1}{3} - \frac{2\, \normalizedmomentcomp{1}^2}{3} & \frac{4\, \normalizedmomentcomp{1}}{3} \end{array}\right)
\end{align}
with eigenvalues
\begin{align*}
\eigenvalue_{1,2} = \frac{2\, \normalizedmomentcomp{1}}{3} \mp \frac{\sqrt{3 - 2\, \normalizedmomentcomp{1}^2}}{3}
\end{align*}
and eigenvectors as in \eqref{eq:MnEigenvectors}. These eigenvalues are shown in \figref{fig:K1Eigenvalues}. Comparing these values with the eigenvalues of the $\MN[1]$ model (dashed), the remarkable difference between minimum-entropy and Kershaw models appears the first time. While the eigenvalues of the $\MN[1]$ model satisfy $\eigenvalue_1(\pm 1) = \eigenvalue_2(\pm 1)$, the $\KN[1]$ eigenvalues evaluate as
\begin{align*}
\eigenvalue_1(1) = \frac13&\neq \eigenvalue_2(1) = 1,\\
\eigenvalue_2(-1) = -\frac13&\neq \eigenvalue_1(-1) = -1.
\end{align*}
A simple calculation shows that
\begin{align*}
\eigenvector_1\cdot \nabla_{\moments} \eigenvalue_1 = -\cfrac{2}{9\momentcomp{0}}\left(2\normalizedmomentcomp{1}+\cfrac{3-\normalizedmomentcomp{1}^2}{\sqrt{3-2\normalizedmomentcomp{1}^2}}\right),
\end{align*}
which is zero for $\normalizedmomentcomp{1} = -1$, but not for $\normalizedmomentcomp{1} = 1$. A similar result is true for $\eigenvector_2\cdot \nabla_{\moments} \eigenvalue_2$, which is zero for $\normalizedmomentcomp{1} = 1$, but not for $\normalizedmomentcomp{1} = -1$. This means that only one of the characteristic fields of the $\KN[1]$ model degenerates on the realizability boundary. In contrast to the $\MN[1]$ model, $\eigenvector_\basisind\cdot\nabla_{\moments} \eigenvalue_\basisind$ is monotonic. The above-derived characteristic fields are shown in \figref{fig:K1Fields}. Note that the author of \cite{Monreal} claimed that at the realizability boundary $\normalizedmomentcomp{1} = \pm 1$ both characteristic fields of the $\KN[1]$ model are linearly degenerate. This is, as shown above, wrong and follows from a mistake in the authors calculation of $\nabla_{\moments} \eigenvalue_i$.
\end{example}
\begin{figure}[htbp]
 \settikzlabel{fig:K1Eigenvalues}\settikzlabel{fig:K1Fields}
\centering
\externaltikz{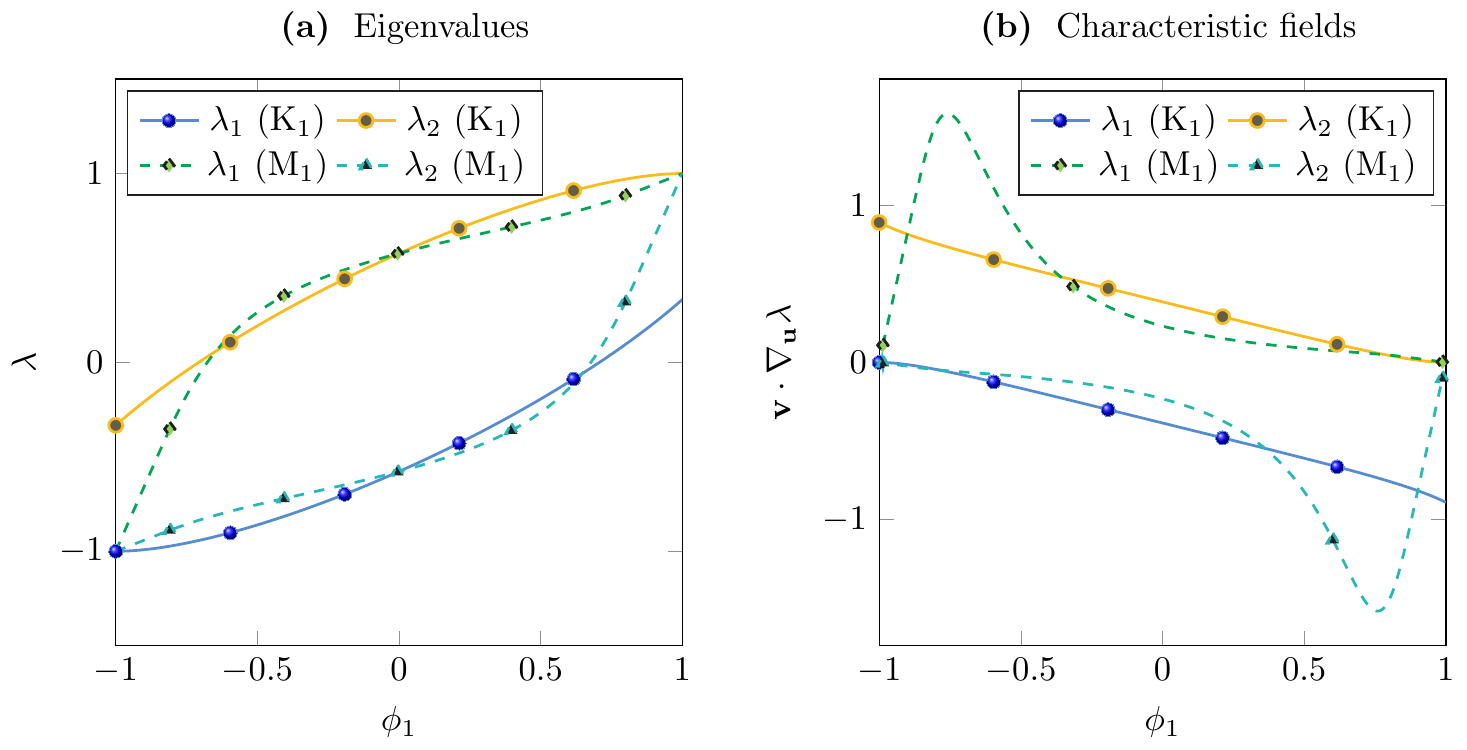}{\relinput{Images/K1EigenvaluesAndFields}}
\caption{{Eigenvalues and characteristic fields ($\momentcomp{0} = 1$) of the $\KN[1]$ (solid) and $\MN[1]$ model (dashed)}.}
\end{figure}

\begin{figure}[htbp]
  \settikzlabel{fig:K2Eigenvalue1}\settikzlabel{fig:K2Eigenvalue2} \settikzlabel{fig:K2Eigenvalue3}
\centering
\externaltikz{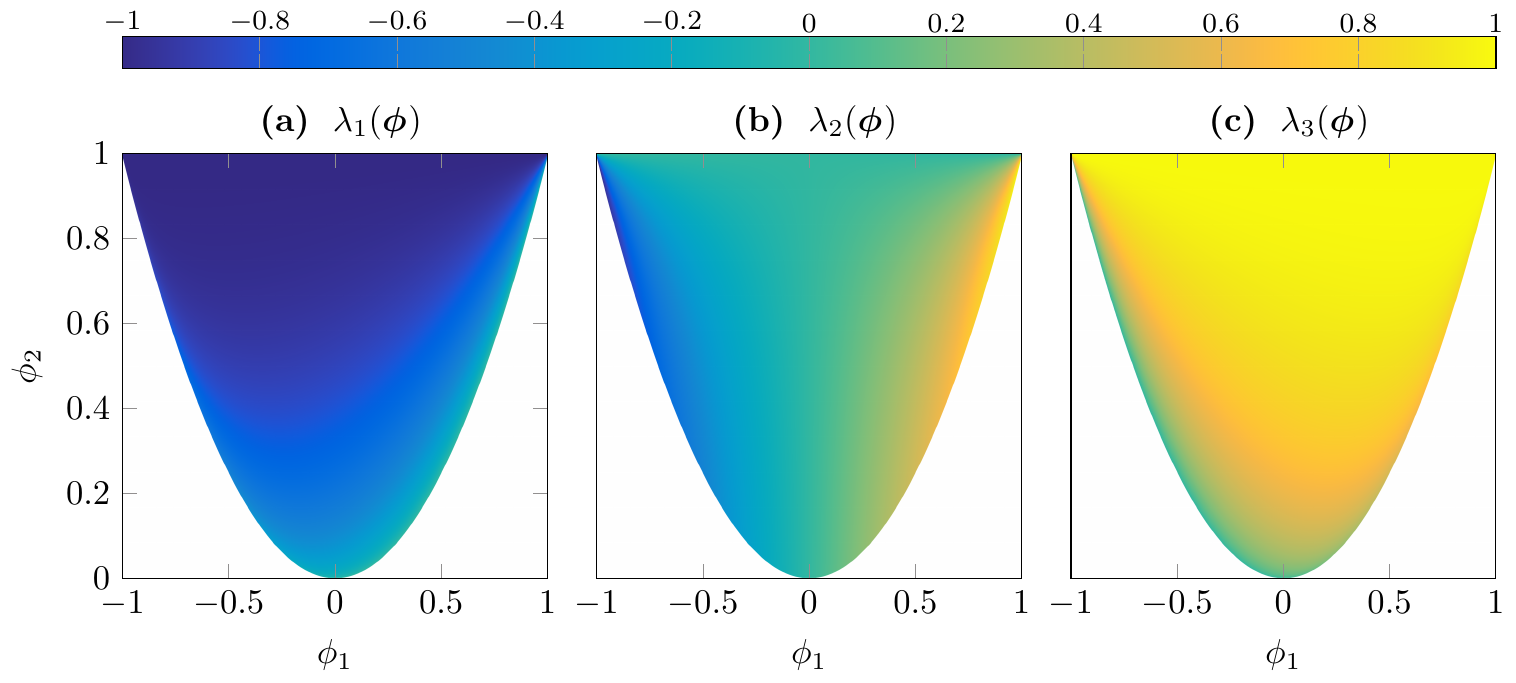}{\relinput{Images/K2Lambdas}}
\caption{Eigenvalues of the flux-Jacobian \eqref{eq:K2Jacobian} for the $\KN[2]$ model.}
\label{fig:K2Eigenvalues}
\end{figure}

\begin{figure}[htbp]
 \settikzlabel{fig:K2MinDiffEigenvalues}\settikzlabel{fig:K2MaxDiffEigenvalues}
\centering
\externaltikz{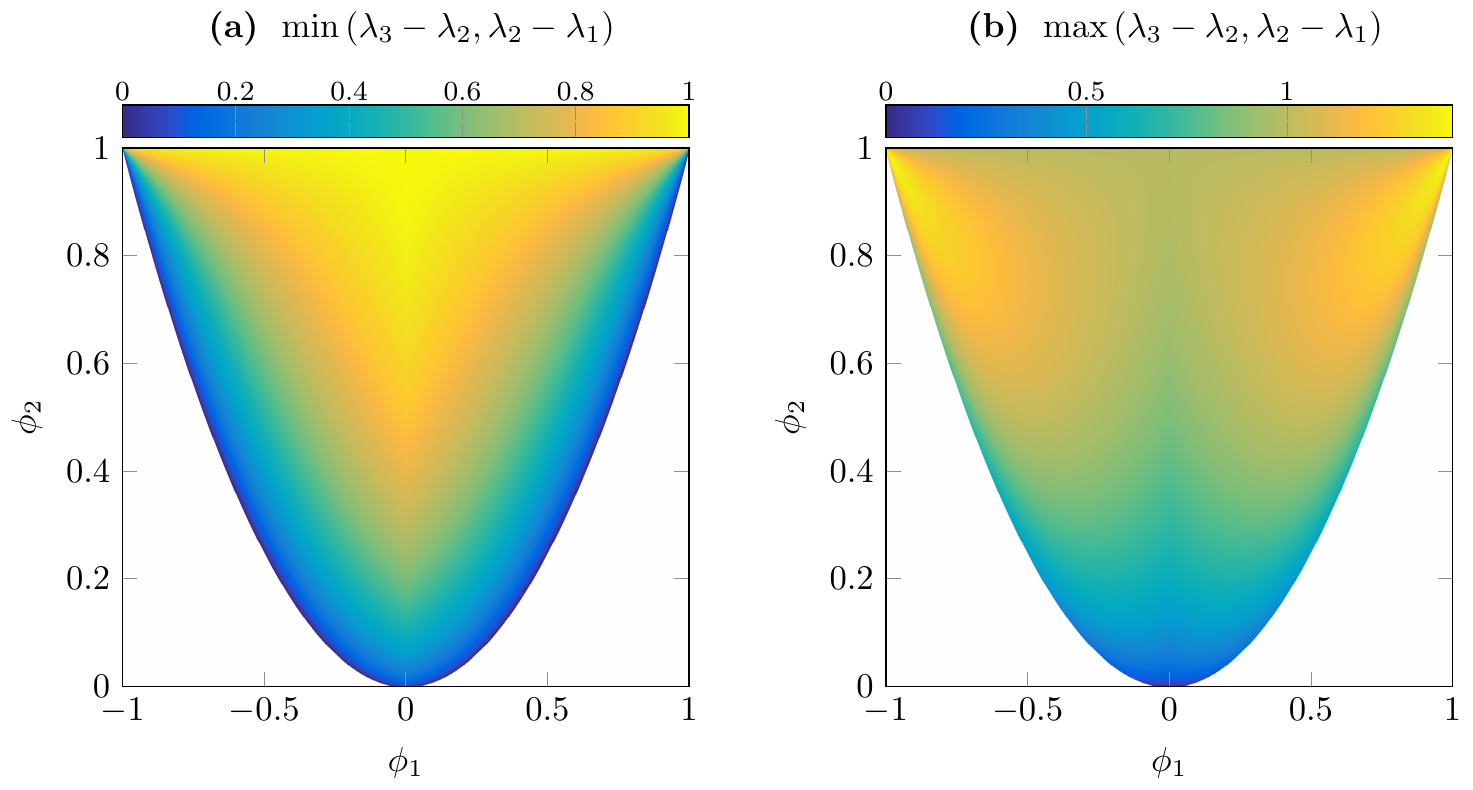}{\relinput{Images/K2Lambdas2}}\\
\caption{Distance of adjacent eigenvalues of the flux-Jacobian \eqref{eq:K2Jacobian} for the $\KN[2]$ model.}
\label{fig:K2Eigenvalues2}
\end{figure}

\begin{figure}[htbp]
 \settikzlabel{fig:K2Field1} \settikzlabel{fig:K2Field2}
\centering
\externaltikz{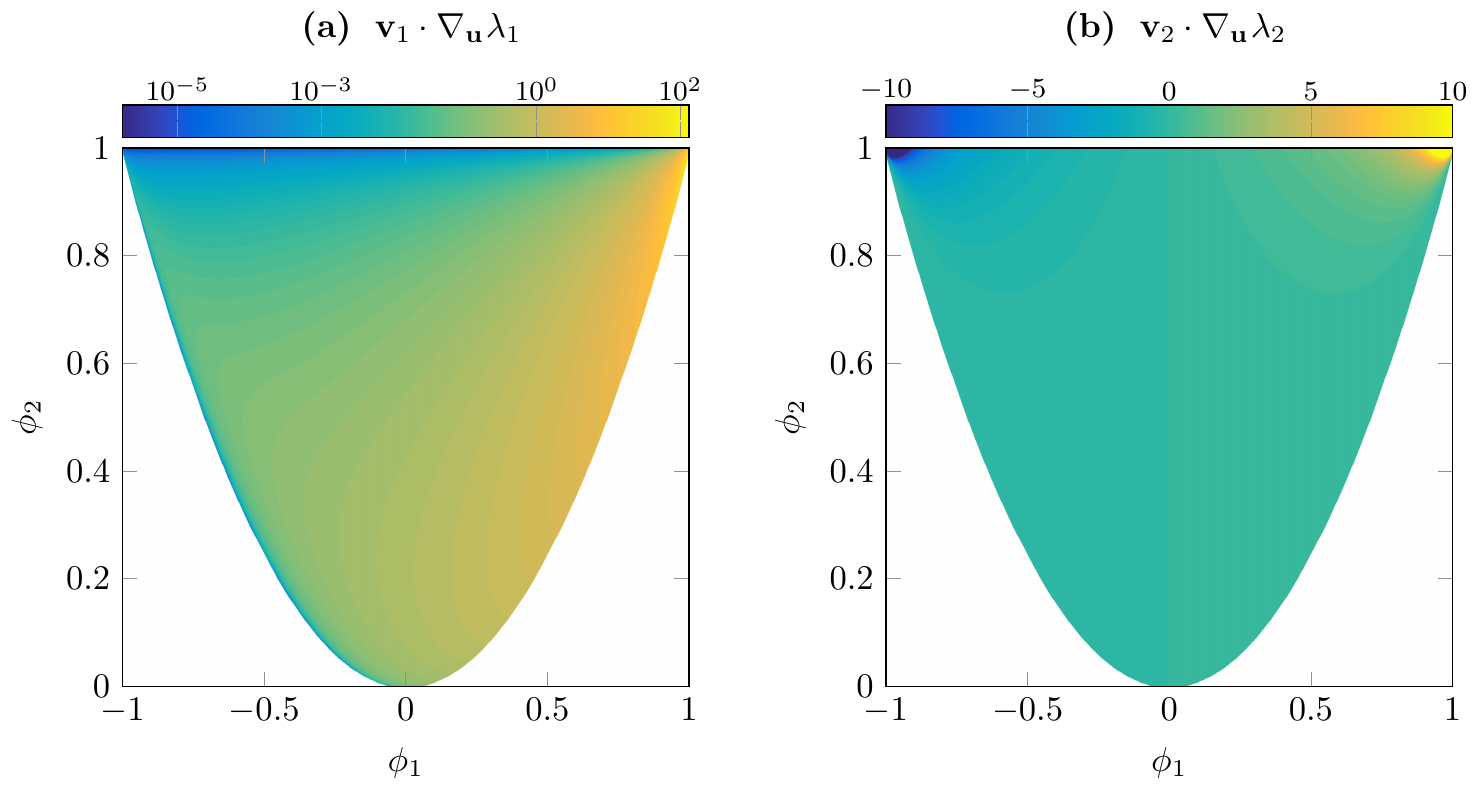}{\relinput{Images/K2Fields}}\\
\caption{Two characteristic fields of the $\KN[2]$ model. The values of the second characteristic field are limited to $[-10,10]$ for easier visualization.}
\label{fig:K2Fields}
\end{figure}

\begin{example}[The {$\KN[2]$} closure]
\label{ex:K2Eigenvalues}
Starting from \eqref{eq:K2closure}, it follows that
\begin{align}
\label{eq:K2Jacobian}
\small\Flux'(\moments) = \left(\begin{array}{ccc} 0 & 1 & 0\\ 0 & 0 & 1\\ \frac{2\, \normalizedmomentcomp{1}\, \left(\normalizedmomentcomp{2} - 1\right)\, \left(\normalizedmomentcomp{2} - \normalizedmomentcomp{1}^2\right)}{{\left(\normalizedmomentcomp{1}^2 - 1\right)}^2} & 1 - \frac{{\left(\normalizedmomentcomp{2} - 1\right)}^2}{2\, {\left(\normalizedmomentcomp{1} + 1\right)}^2} - \frac{{\left(\normalizedmomentcomp{2} - 1\right)}^2}{2\, {\left(\normalizedmomentcomp{1} - 1\right)}^2} & \frac{2\, \normalizedmomentcomp{1}\, \left(\normalizedmomentcomp{2} - 1\right)}{\normalizedmomentcomp{1}^2 - 1} \end{array}\right).
\end{align}
Since the general formula for the eigenvalues is very lengthy, it is omitted here. Instead, the eigenvalues are plotted in \figref{fig:K2Eigenvalues}. 
%The zeros of the eigenvalues are given by
%\begin{align*}
%\lambda_1(\normalizedmoments) &= 0 \Leftrightarrow \normalizedmomentcomp{1}\in \left[0,1\right), \normalizedmomentcomp{2} = \normalizedmomentcomp{1}^2\\
%\lambda_2(\normalizedmoments) &= 0 \Leftrightarrow \normalizedmomentcomp{1}=0 \text{ or } \normalizedmomentcomp{2} = 1\\
%\lambda_3(\normalizedmoments) &= 0 \Leftrightarrow \normalizedmomentcomp{1}\in \left(-1,0\right], \normalizedmomentcomp{2} = \normalizedmomentcomp{1}^2
%\end{align*}

The minimal and maximal distance between two adjacent eigenvalues is depicted in \figref{fig:K2Eigenvalues2}. It is visible that on the lower realizability boundary, where $\normalizedmomentcomp{2} = \normalizedmomentcomp{1}^2$ and
\begin{align}
\label{eq:K2JacobianOnLowerBoundary}
\Flux'(\moments) = \left(\begin{array}{ccc} 0 & 1 & 0\\ 0 & 0 & 1\\ 0 & - \normalizedmomentcomp{1}^2 & 2\, \normalizedmomentcomp{1} \end{array}\right),
\end{align}
the minimal distance is zero, while the maximal distance is greater than zero almost everywhere. Calculating the eigenvalues of \eqref{eq:K2JacobianOnLowerBoundary}, one obtains zero as a single eigenvalue and $\normalizedmomentcomp{1}$ twice. For $\normalizedmomentcomp{1}=\normalizedmomentcomp{2}=0$, all eigenvalues coincide.

On the other hand, looking at the upper realizability boundary, i.e. $\normalizedmomentcomp{2}=1$, it follows that
\begin{align}
\label{eq:K2JacobianOnUpperBoundary}
\Flux'(\moments) = \left(\begin{array}{ccc} 0 & 1 & 0\\ 0 & 0 & 1\\ 0 & 1 & 0 \end{array}\right),
\end{align}
which has eigenvalues $\eigenvalue_1 = -1$, $\eigenvalue_2 = 0$ and $\eigenvalue_3 = 1$. Therefore, the system is, as the $\KN[1]$ model, strictly hyperbolic on the upper realizability boundary. Note that this calculation also implies that the eigenvalues of \eqref{eq:K2Jacobian} are discontinuous at $\abs{\normalizedmomentcomp{1}}=\normalizedmomentcomp{2}=1$. This is different for the $\KN[2]$ model derived in \cite{Monreal}. However, the model there is linearly degenerate everywhere. Nevertheless, the $\MN[2]$ model shows the same behaviour, see e.g. \cite{Wright2009} where this is shown for a different entropy.
\end{example}
The characteristic fields, given in \figref{fig:K2Fields}, are always zero on exactly two parts of the realizability boundary. For example, $\eigenvector_1\cdot\nabla_{\moments}\eigenvalue_1$ is zero for $\normalizedmomentcomp{2}=\normalizedmomentcomp{1}^2$ and $\normalizedmomentcomp{1}\in[-1,0]$ or $\normalizedmomentcomp{2}=1$. Similarly, $\eigenvector_3\cdot\nabla_{\moments}\eigenvalue_3 = 0$ for $\normalizedmomentcomp{2}=\normalizedmomentcomp{1}^2$ and $\normalizedmomentcomp{1}\in[0,1]$ or $\normalizedmomentcomp{2}=1$. The second field (\figref{fig:K2Field2}) is zero for $\normalizedmomentcomp{2}=\normalizedmomentcomp{1}^2$ or $\normalizedmomentcomp{1} = 0$. 
Overall, at least one of the characteristic fields is never linearly degenerate, which is in good coincidence with the behaviour of the $\KN[1]$ system (see \figref{fig:K1Fields}). Note that the third characteristic field is skew-symmetric to the first field along the symmetry axis $\normalizedmomentcomp{1} = 0$.

\section{Numerical experiments}
This section contains some often-used benchmark problems for moment models. With them it is possible to qualitatively investigate the differences between the derived moment models. The following results are calculated on a highly-resolved first-order grid ($\ncells = 10000$), using the techniques given in \cite{Schneider2015a}. Everything is calculated with isotropic scattering, i.e. $\collision{\distribution} = \frac12\momentcomp{0}-\distribution$. The \emph{reference solution} is given by the $\PN[199]$ model.

\subsection{Plane source}
\label{sec:Planesource}
In this test case an isotropic distribution with all mass concentrated in the middle of an infinite domain $\z \in
(-\infty, \infty)$ is defined as initial condition, i.e.
\begin{align*}
 \distributiontzero(\z, \SCheight) = \distributionvacuum + \delta(\z),
\end{align*}
where the small parameter $\distributionvacuum = 0.5 \times 10^{-8}$ is used to
approximate a vacuum.
In practice, a bounded domain must be used which is large
enough that the boundary should have only negligible effects on the
solution. For the final time $\tf = 1$, the domain is set to $\Domain = [-1.2, 1.2]$ (recall that for all presented models the maximal speed of propagation is bounded in absolute value by one).

At the boundary the vacuum approximation
\begin{align*}
 \distributionboundary(\timevar,\zL,\SCheight) \equiv \distributionvacuum \quand
 \distributionboundary(\timevar,\zR,\SCheight) \equiv \distributionvacuum
\end{align*}
 is used again. Furthermore, the physical coefficients are set to $\scattering \equiv 1$, $\absorption \equiv 0$ and $\source \equiv 0$.

All solutions are computed with an even number of cells, so the initial Dirac delta lies on a cell boundary.
Therefore it is approximated by splitting it into the cells immediately to the left and right. In all figures below, only positive $\z$ are shown since the solutions are always symmetric around $\z = 0$.\\

Noting that the method of moments is indeed a type of spectral method, it can be expected that due to the non-smoothness of the initial condition the convergence towards the kinetic solution of this test case is slow (note that $\distributiontzero(\cdot, \SCheight)\notin\Lp{p}$ for any $p$).

\begin{figure}[htbp]
\centering
\externaltikz{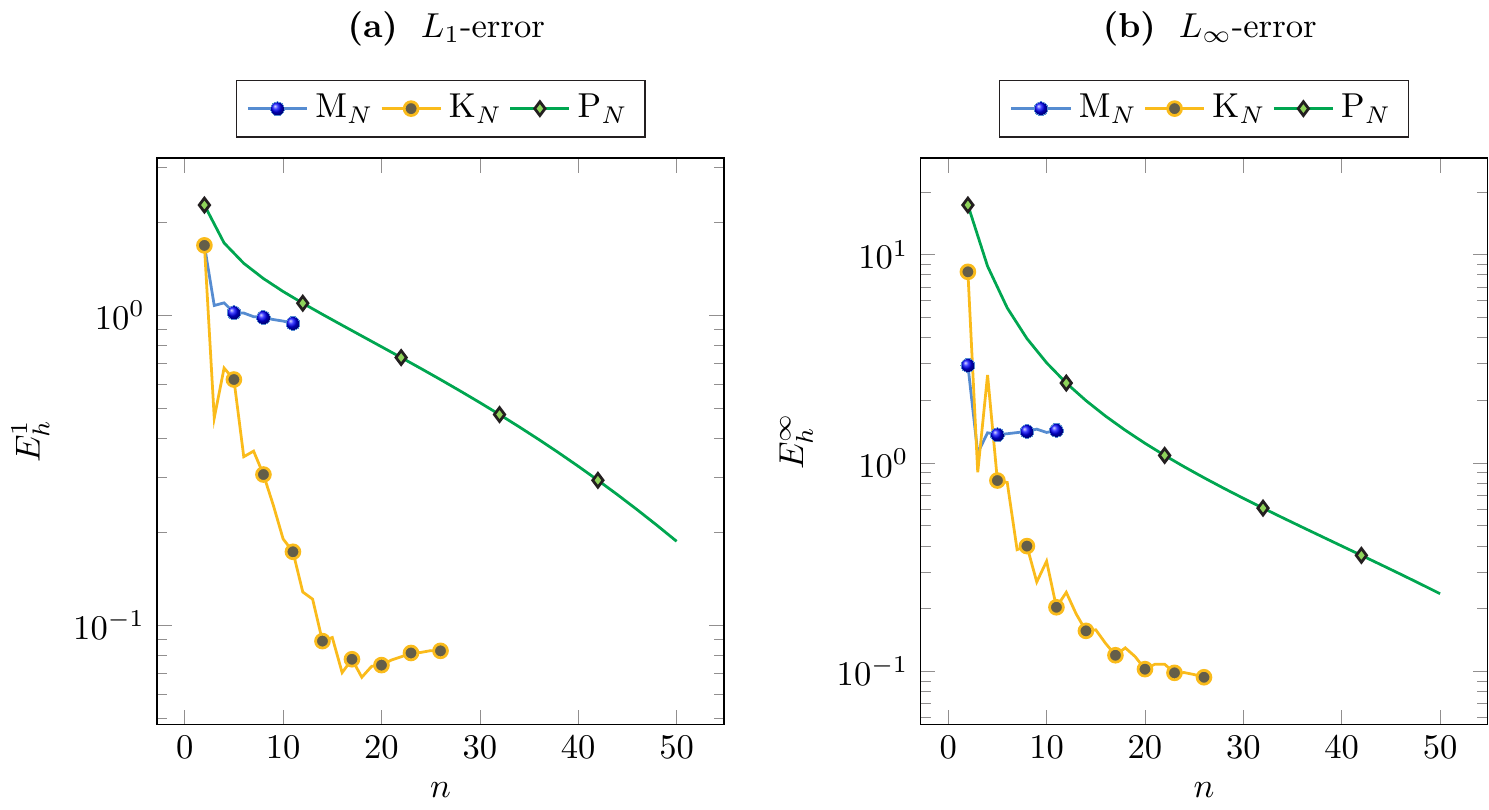}{\relinput{Images/PlanesourceConvergence}}
\caption{{$\Lp{p}$-convergence of $\momentcomp{0}$ for different models against the $\PN[199]$ reference solution of the plane-source test with isotropic scattering depending on the number of moments $\momentnumber = \momentorder+1$. Errors are evaluated at the final time $\tf = 1$.}}
\label{fig:PlanesourceConvergence}
\end{figure}

This is shown in \figref{fig:PlanesourceConvergence}, where the $\Lp{1}$- and $\Lp{\infty}$-error measured in the local particle density evaluated at the final time are plotted against the number of moments for the isotropic-scattering operator. Minimum-entropy models behave reasonably better than the standard $\PN$ method. This has been observed as well in \cite{Hauck2010}. 

A big surprise is the behaviour of the Kershaw models. Although there is a structural similarity between them and their corresponding minimum-entropy versions, the difference obviously increases with higher moment order. 

\begin{figure}[htbp]
\centering
\settikzlabel{fig:PlanesourceIsotropicCutsMn}\settikzlabel{fig:PlanesourceIsotropicCutsKn}
\externaltikz{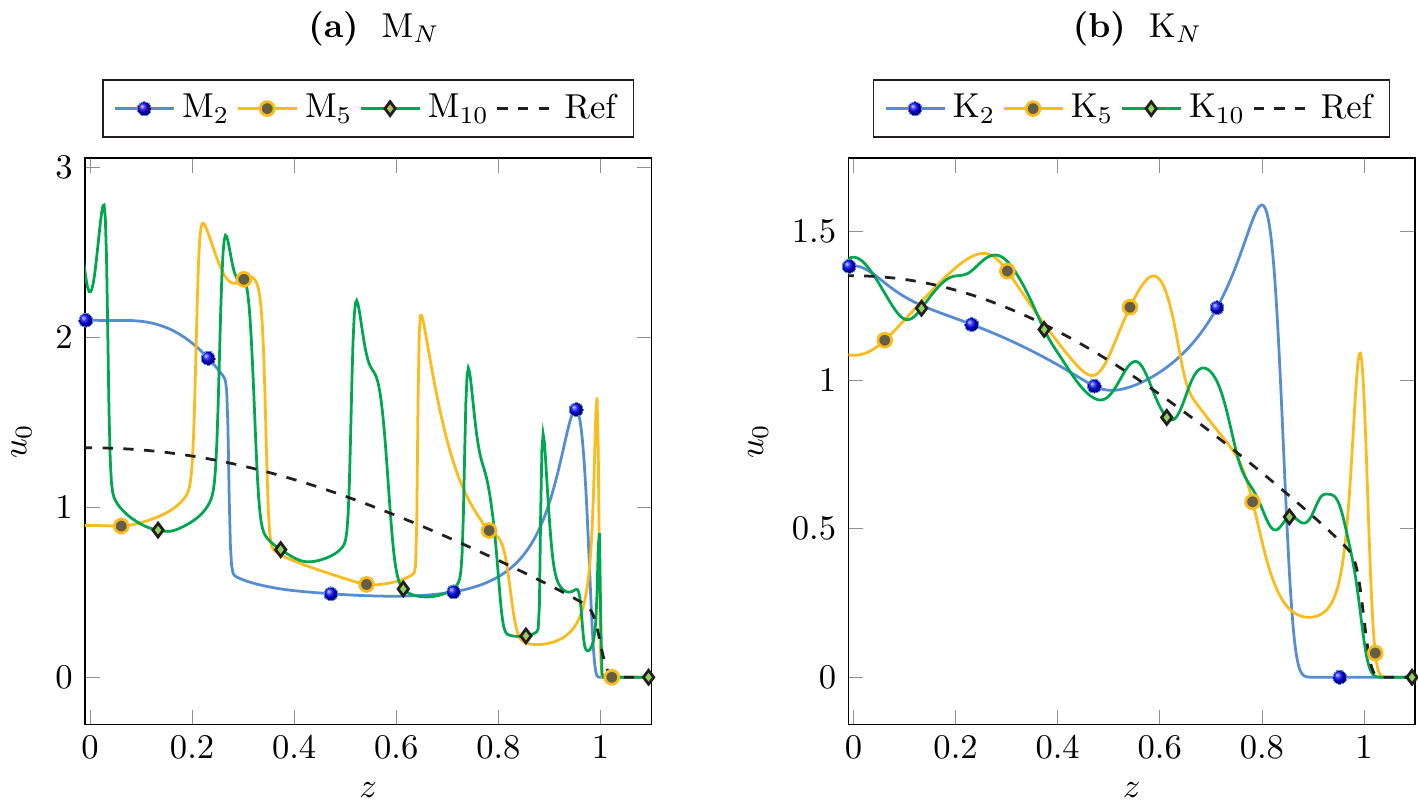}{\relinput{Images/PlanesourceIsotropicCutsFullMoments}}
\caption{{$\MN$ and $\KN$ solutions of the plane-source test at $\tf = 1$ with isotropic scattering.}}
\label{fig:PlanesourceIsotropicCutsFullMoments}
\end{figure}

This can be seen in \figref{fig:PlanesourceIsotropicCutsFullMoments} where $\MN$ (\figref{fig:PlanesourceIsotropicCutsMn}) and $\KN$ models (\figref{fig:PlanesourceIsotropicCutsKn}) of multiple orders are plotted beside each other.

\subsection{Source beam}
\label{sec:SourceBeam}
Finally, a discontinuous version of the source-beam problem from
\cite{Hauck2013} is presented.
The spatial domain is $\Domain = [0,3]$, and
\begin{gather*}
 \absorption(\z) = \begin{cases}
   1 & \text{ if } \z\leq 2,\\
   0 & \text{ else},
  \end{cases} \quad
 \scattering(\z) = \begin{cases}
   0 & \text{ if } \z\leq 1,\\
   2 & \text{ if } 1<\z\leq 2,\\
   10 & \text{ else}
  \end{cases} \quad
 \source(\z) = \begin{cases}
   1 & \text{ if } 1\leq \z\leq 1.5,\\
   0 & \text{ else},
  \end{cases}
\end{gather*}
with initial and boundary conditions
\begin{gather*}
 \distributiontzero(\z, \SCheight) \equiv \distributionvacuum, \\
 \distributionboundary(\timevar,\zL,\SCheight) = \cfrac{e^{-10^5(\SCheight-1)^2}}{\ints{e^{-10^5(\SCheight-1)^2}}}
 \quand
 \distributionboundary(\timevar,\zR,\SCheight) \equiv \distributionvacuum.
\end{gather*}
The final time is $\tf = 2.5$ and the same vacuum approximation $\distributionvacuum$ as in the plane-source problem is used.

\begin{figure}[htbp]
\centering
\externaltikz{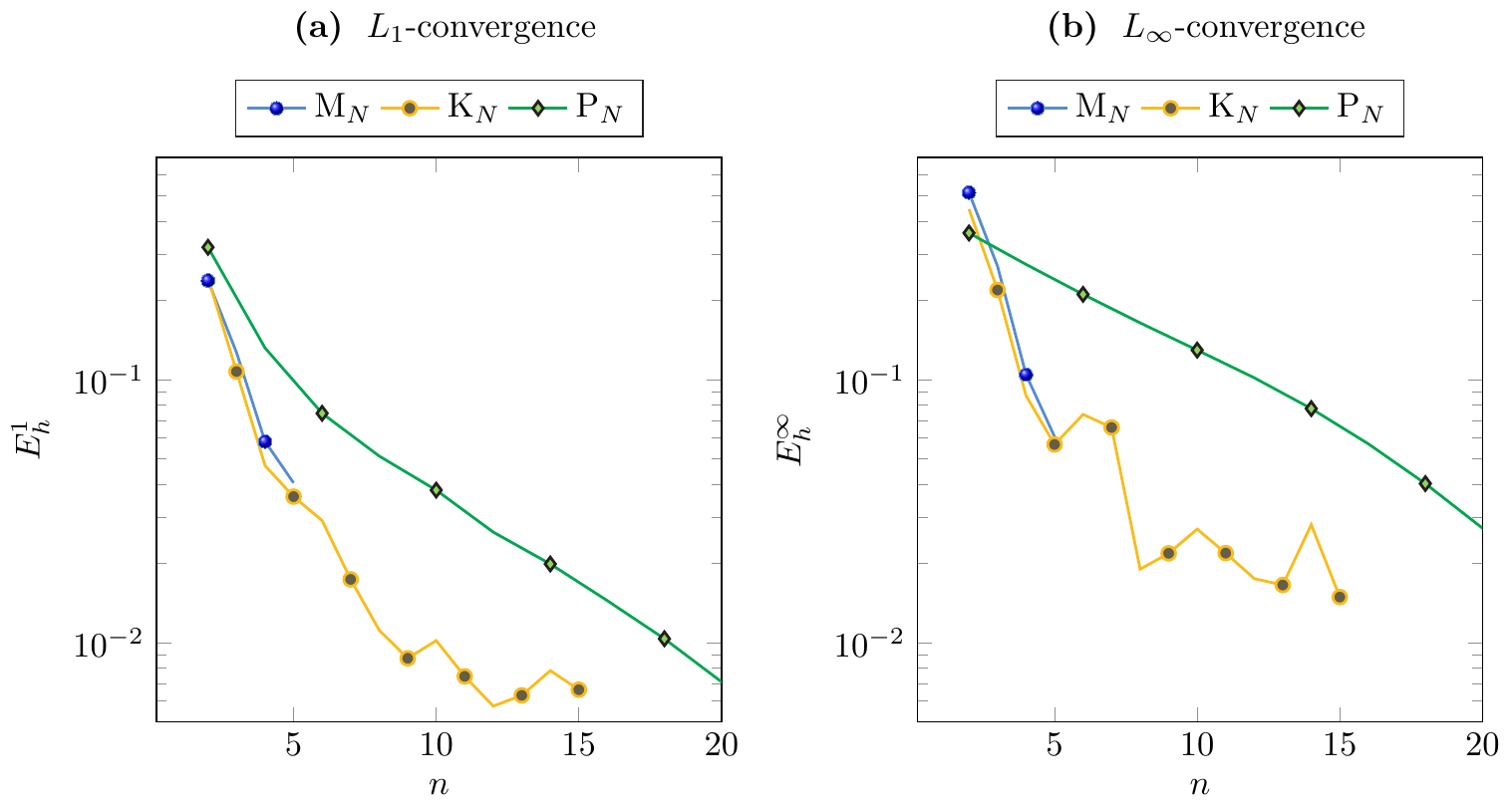}{\relinput{Images/SourceBeamConvergence}}
\caption{$\Lp{p}$-convergence of $\momentcomp{0}$ for different models against the $\PN[199]$ reference solution of the source-beam test case with isotropic scattering depending on the number of moments $\momentnumber = \momentorder+1$. Errors are evaluated at the final time $\tf = 2.5$.}
\label{fig:SourceBeamConvergence}
\end{figure}

$\Lp{p}$-errors are shown in \figref{fig:SourceBeamConvergence}. It is visible that $\MN$ and $\KN$ models perform similarly. Both models outperform the classical $\PN$ model.

\begin{figure}[h!]
\centering
\settikzlabel{fig:SourceBeamIsotropicCutsMn}\settikzlabel{fig:SourceBeamIsotropicCutsKn}
\externaltikz{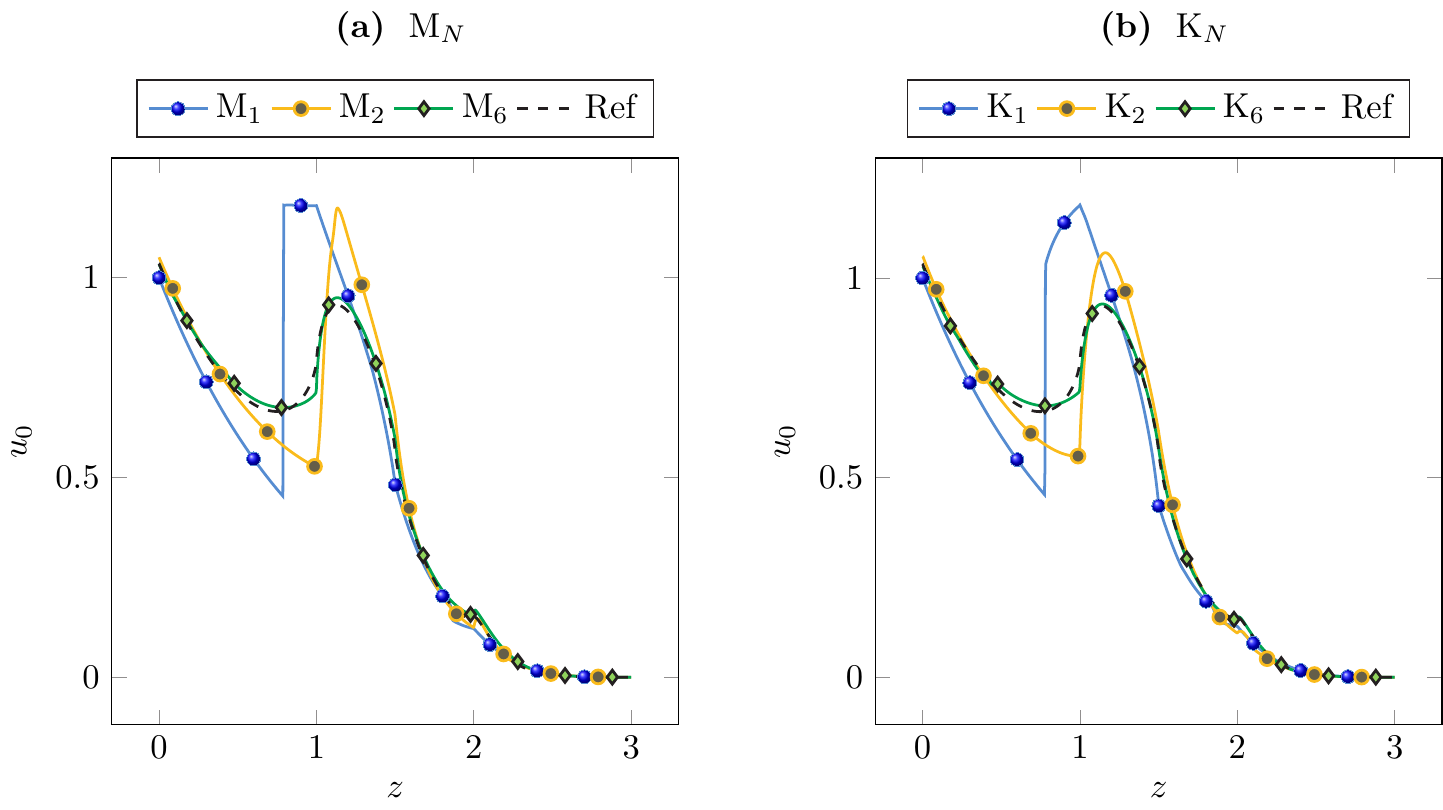}{\relinput{Images/SourceBeamIsotropicCutsFullMoments}}
\caption{$\MN$ and $\KN$ solutions of the source-beam test at $\tf = 2.5$ with isotropic scattering.}
\label{fig:SourceBeamIsotropicCutsFullMoments}
\end{figure}

\section{Conclusions and outlook}
In this paper the basic concepts of deriving full-moment Kershaw closures were derived. These models provide a huge gain in efficiency compared to the state-of-the-art minimum-entropy models, since they can be closed (in principle) analytically using the available realizability theory. Benchmark tests confirm that Kershaw closures can indeed compete with minimum-entropy models. In some situations, they are even better.

Although the gain in efficiency is very big, recent results for minimum-entropy models showed that using high-order numerics is still advantageous \cite{Schneider2015a,Schneider2015b}. The big problem in high-order schemes for moment models is that the property of realizability has to be preserved during the simulation, since otherwise the fluxes cannot be evaluated. Future work will have to investigate how to adapt the scheme in \cite{Schneider2015a} to Kershaw closures. Furthermore, different scattering operators should be taken into account, like the slightly more complicated (in terms of realizability preservation) Laplace-Beltrami operator.

Finally, the concepts have to be lifted to higher dimensions. While fully three-dimensional first-order variants of Kershaw closures exist \cite{Ker76,Schneider2015c}, no higher-order models or a completely closed theory is available.  

% Bibliography
%%%%%%%%%%%%%%
\bibliographystyle{siam}
\bibliography{bibliography}

\end{document}